\documentclass[11pt]{article}

\usepackage{fullpage}
\usepackage{authblk}
\usepackage{natbib}
\usepackage{algorithm}
\usepackage{algpseudocode}
\usepackage{amsmath,amsthm,amsfonts}
\usepackage{amssymb}
\usepackage[subnum]{cases}
\usepackage{tcolorbox}
\usepackage{mathrsfs}
\usepackage{amssymb}
\usepackage{color}
\usepackage{mathrsfs}
\usepackage{enumitem}
\usepackage{bm}
\usepackage{multirow}
\usepackage{booktabs}
\usepackage{makecell}
\usepackage{graphicx}
\usepackage{epstopdf}
\usepackage{subcaption}
\usepackage{comment}
\usepackage{cases}
\usepackage{appendix}
\usepackage{tikz}
\usetikzlibrary{arrows,shapes}
\usepackage[colorlinks= true, linkcolor=brown, citecolor=blue, urlcolor=black]{hyperref}
\usepackage{cleveref}
\usepackage{marginnote}
\usepackage{diagbox} 
\usepackage{multirow} 
\usepackage{nicematrix,tikz}
\usetikzlibrary{fit}

\usepackage[utf8]{inputenc}
\usepackage{amsmath}
\usepackage{nicematrix}
\usepackage{bm}
\usepackage{blkarray}

\usepackage{arydshln}
\usepackage{xurl}   
\usepackage{seqsplit} 

\numberwithin{equation}{section}

\theoremstyle{definition}

\newtheorem{theorem}{Theorem}[section]
\newtheorem{lemma}[theorem]{Lemma}

\newtheorem{defn}[theorem]{Definition}

\usepackage{booktabs}
\usepackage{tikz}
\usepackage{pgfplots}
\usetikzlibrary{calc,intersections}

\usetikzlibrary{shapes.geometric, arrows, positioning}

\tikzset{
	mybox/.style  = {draw, rectangle, minimum width=4cm, minimum height=0.8cm, text centered, text width=4.4cm,   
		font=\normalsize},
	box/.style  = {draw, rectangle, minimum width=2.0cm, minimum height=0.6cm, text centered, text width=3.0cm,   
		font=\normalsize},
	myarrow/.style = {line width=0.2pt, draw=black, -triangle 60, postaction={draw, line width=0.2pt, shorten >=10pt,-}}
}

\tikzstyle{arrow} = [->, >=stealth, -triangle 60]

\allowdisplaybreaks

\makeatletter
\newcommand{\leqnomode}{\tagsleft@true}
\newcommand{\reqnomode}{\tagsleft@false}
\makeatother

\begin{document}

\title{On Pseudospectral Concentration for Rank-$1$ Sampling}

\author[1]{Kuo Gai}
\author[2,3]{Bin Shi\thanks{Corresponding author, Email: \url{binshi@fudan.edu.cn}}} 
%
\affil[1]{School of Artificial Intelligence, Wuhan University, Wuhan 430072, China}
\affil[2]{Center for Mathematics and Interdisciplinary Sciences, Fudan University, Shanghai 200433, China}
\affil[3]{Shanghai Institute for Mathematics and Interdisciplinary Sciences, Shanghai 200438, China}

\date\today

\maketitle

%

\begin{abstract}
Pseudospectral analysis serves as a powerful tool in matrix computation and the study of both linear and nonlinear dynamical systems. Among various numerical strategies, random sampling, especially in the form of rank-$1$ perturbations, offers a practical and computationally efficient approach.  Moreover, due to invariance under unitary similarity, any complex matrix can be reduced to its upper triangular form, thereby simplifying the analysis. In this study. we develop a quantitative concentration theory for the pseudospectra of complex matrices under rank-$1$ random sampling perturbations, establishing a rigorous probabilistic framework for spectral characterization. First, for normal matrices, we derive a regular concentration inequality and demonstrate that the separation radius scales with the dimension as $\delta_d \sim 1/\sqrt{d}$. Next, for the equivalence class of nilpotent Jordan blocks,  we exploit classical probabilistic tools, specifically, the Hanson–Wright concentration inequality and the Carbery–Wright anti-concentration inequality, to obtain singular concentration bounds, and demonstrate that the separation radius exhibits the same dimension-dependent scaling. This yields a singular pseudospectral concentration framework. Finally, observing that upper triangular Toeplitz matrices can be represented via the symbolic polynomials of nilpotent Jordan blocks, we employ partial fraction decomposition of rational functions to extend the singular framework to the equivalence class of upper triangular Toeplitz matrices.
\end{abstract}

\section{Introduction}
\label{sec: intro}

Eigenvalues, and more broadly, the spectrum of a matrix or linear operator, have long played a central role in analyzing both linear and nonlinear systems across science and engineering~\citep{golub2013matrix}. However, for non-Hermitian, or more generally, nonnormal matrices, where eigenvectors are not orthogonal, the spectrum alone often fails to yield meaningful insights. This limitation becomes especially pronounced in high-dimensional or strongly nonnormal systems, where even small perturbations can trigger large and unpredictable responses. Such behavior arises in diverse fields, including hydrodynamics~\citep{reddy1993energy, reddy1993pseudospectra, trefethen1993hydrodynamic}, plasma physics~\citep{borba1994pseudospectrum}, atmospheric and oceanic dynamics~\citep{farrell1993stochastic, farrell1992adjoint}, control theory~\citep{hinrichsen1994stability}, and numerical algorithms~\citep{van1993linear}. To better characterize such sensitivity, the concept of \textit{pseudospectra} has emerged as a powerful extension of classical eigenvalue analysis~\citep{trefethen1997pseudospectra, trefethen2020spectra}.  Pseudospectral analysis examines how the spectrum responds to small perturbations, identifying regions in the complex plane where the resolvent norm becomes large.

Despite its conceptual appeal, computing pseudospectra remains challenging~\citep{trefethen1999computation}.  Traditional grid-based methods involve evaluating the resolvent norm over a dense mesh in the complex plane, which requires the solution of many shifted linear systems~\citep[Section 39]{trefethen2020spectra}. While straightforward in principle, these methods become computationally expensive and often impractical for large-scale or highly nonnormal matrices. As a practical alternative,  sampling-based techniques, referred to as~``\textit{Poor Man's Algorithms}'' by~\citet[Section 39]{trefethen2020spectra}, stand out for their simplicity and efficiency. By randomly sampling perturbations, these methods empirically approximate the pseudospectrum, identifying regions of potential spectral instability while avoiding the excessive computational overhead of grid-based approaches.



From a computational standpoint, we are guided by the second equivalent definition of pseudospectra presented in~\citet[Section 2]{trefethen2020spectra}; see~\Cref{defn: pseudospectra} for the precise statement. In practice, one can generate normalized full-rank perturbations via random sampling, as demonstrated in~\citet{trefethen1999computation}):\footnote{\noindent Throughout this paper, unless otherwise stated, all vectors $\pmb{\xi}$ are in the $d$-dimensional complex space $\mathbb{C}^d$, and all matrices $E$ are $d\times d$ complex-valued matrices in $\mathbb{C}^{d \times d}$. The norm $\|\cdot\|$ refers to the spectral norm (also known as the 2-norm), For a vector $\pmb{\xi} \in \mathbb{C}^d$, we define $\|\pmb{\xi}\| = \sqrt{\pmb{\xi}^{\dagger} \pmb{\xi}}$, where $\dagger$ denotes the conjugate transpose. For a matrix $E \in \mathbb{C}^{d \times d}$, the norm is $\|E\| = \max \sqrt{\lambda(E^{\top}E)}$, which is the largest singular value of $E$.}  
\begin{tcolorbox}
\begin{itemize}
\item[1.] \textbf{Full-Rank Sampling:} Generate a matrix  $E \in \mathbb{C}^{d \times d}$ with i.i.d. entries from the standard complex normal distribution $\mathcal{CN}(0,1)$, and then normalize it to unit spectral norm, i.e., $E/\|E\| $.
\end{itemize}
\end{tcolorbox}
\noindent As established in~\citet[Theorem 2.1]{trefethen2020spectra}, the definition of pseudospectra remains valid when full-rank perturbations are replaced with rank-$1$ perturbations. A rigorous justification is provided in~\Cref{thm: rank1-pseudospectra}. This observation motivates the use of a more  computationally efficient rank-1 sampling method:
\begin{tcolorbox}
\begin{itemize}
\item[2.] \textbf{Rank-1 Sampling:} Independently draw two complex normal vectors $\pmb{u}, \pmb{v} \sim \mathcal{CN}(\pmb{0}, I)$, and construct the rank-$1$ matrix as
\[
E = \frac{ \pmb{u}  \pmb{v}^{\dagger}}{\|\pmb{u}\| \|\pmb{v}\|},
\]
which ensures $\|E\| = 1$ while  offering significant computational advantages.
\end{itemize}
\end{tcolorbox}
\noindent The computational efficiency of rank-1 sampling, first highlighted by~\citep{riedel1994generalized}, provides substantial speedups over full-rank sampling. Normalizing a full-rank matrix requires $O(d^3)$ operations, whereas the rank-1 approach reduces this to $O(d^2)$, as emphasized in~\citet[Section 39]{trefethen2020spectra}. Further support comes from~\citet{trefethen1999computation}, which demonstrated the effectiveness of rank-1 perturbations in capturing the pseudospectral structure of a matrix.
\subsection{Random perturbations: invariance under unitary similarity}
\label{subsec: unitary}

A key observation in pseudospectral analysis is that random perturbation matrices, whether generated via  full-rank or rank-1 sampling, are invariant under unitary similarity transformations. This invariance plays a crucial role in understanding the structural behavior of pseudospectra under random perturbations. Recall two matrices, $E, E_1 \in \mathbb{C}^{d \times d}$ are said to be unitarily similar if
\[
E_1 = U^{\dagger} E U, 
\]
where $U \in \mathbb{C}^{d \times d}$ is a unitary matrix. If $E$ is sampled with i.i.d. entries from the standard complex normal distribution $\mathcal{CN}(0,1)$, the rotational invariance of this distribution ensures that $E_1$ inherits the same statistical properties. Moreover, the spectral norm is preserved under unitary similarity: 
\[
\|E_1\| = \max \sqrt{\lambda(E_1^{\dagger}E_1)} = \max \sqrt{\lambda(U^{\dagger} E^{\dagger} E U)} = \max \sqrt{\lambda(E^{\dagger} E)} = \|E\|, 
\]
where the third equality follows from the unitary invariance of eigenvalues. Hence, normalized full-rank perturbations are invariant under unitary similarity:
\begin{equation}
\label{eqn: unitary-similar-full-rank}
\frac{E_1}{\|E_1\|} = U^{\dagger} \left(\frac{E}{\|E\|}\right) U. 
\end{equation}
A similar conclusion holds for rank-1 sampling. Let $\pmb{u}, \pmb{v} \sim \mathcal{CN}(\pmb{0}, \pmb{I})$ be two independent complex normal vectors in $\mathbb{C}^2$. Since unitary transformations preserve both the distribution and norm (i.e., $U^{\dagger} \pmb{u} \sim \mathcal{CN}(\pmb{0}, \pmb{I}_{d \times d})$ and $\| U^{\dagger} \pmb{u}\| = \| \pmb{u}\|$), the resulting rank-1 perturbation matrix is also invariant under unitary similarity as
\begin{equation}
\label{eqn: unitary-similar-rank-1}
E_1 = U^{\dagger} E U =  \frac{ \left(U^{\dagger} \pmb{u}\right)  \left( U^{\dagger} \pmb{v}\right)^{\dagger}}{\|\pmb{u}\| \|\pmb{v}\|} = \frac{ \left(U^{\dagger} \pmb{u}\right)  \left( U^{\dagger} \pmb{v}\right)^{\dagger}}{\|U^{\dagger} \pmb{u}\| \| U^{\dagger} \pmb{v}\|}.
\end{equation}

This invariance highlights the robustness of both sampling strategies, full-rank and rank-1, to changes in the underlying unitary basis, reinforcing their practical utility in pseudospectral analysis.  By exploiting this property, we can simplify the analysis of random perturbations through unitary equivalence classes. In particular, since both sampling methods are invariant under unitary similarity (as shown in~\eqref{eqn: unitary-similar-full-rank} and~\eqref{eqn: unitary-similar-rank-1}), we can invoke the Schur decomposition (\citep[Theorem 7.1.3]{golub2013matrix}),  which states that every matrix is unitarily similar to an upper triangular matrix. This insight allows us to reduce the computation of pseudospectra to a canonical form. Specifically, for any $\epsilon > 0$,   there exists a unitary matrix $U \in \mathbb{C}^{d \times d}$ such that the perturbed matrix satisfies
\begin{equation}
\label{eqn: perturb-unitary}
U^{\dagger}(A+\epsilon E)U = T + \epsilon U^{\dagger}EU = T + \epsilon E_1,
\end{equation}
where $T$ is an upper triangular matrix, and $E_1$ inherits the same statistical properties as $E$. Thus, to understand the pseudospectral behavior of arbitrary matrices under random perturbations, it suffices to analyze the case of upper triangular matrices.  A rigorous justification for this reduction is presented in~\Cref{sec: prelim}. 


\subsection{Pseudospectral phenomena in rank-1 sampling }
\label{subsec: pseudospectral-phenomenon}

For the remainder of this paper, we focus exclusively on rank-1 sampling, motivated by its substantial computational advantages over full-rank sampling. As demonstrated in the pioneering work of~\citet{trefethen1999computation} and further formalized in~\citet{trefethen2020spectra}, rank-1 sampling provides an efficient and reliable methods of probing pseudospectral behavior.  Consequently, we restrict our attention to rank-1 sampling perturbations and omit further discussion of the full-rank case.

\subsubsection{Regular pseudospectral concentration}
\label{subsubsec: regular-concentration}

\begin{figure}[htb!]
\centering
\begin{subfigure}[t]{0.45\linewidth}
\centering
\includegraphics[scale=0.48]{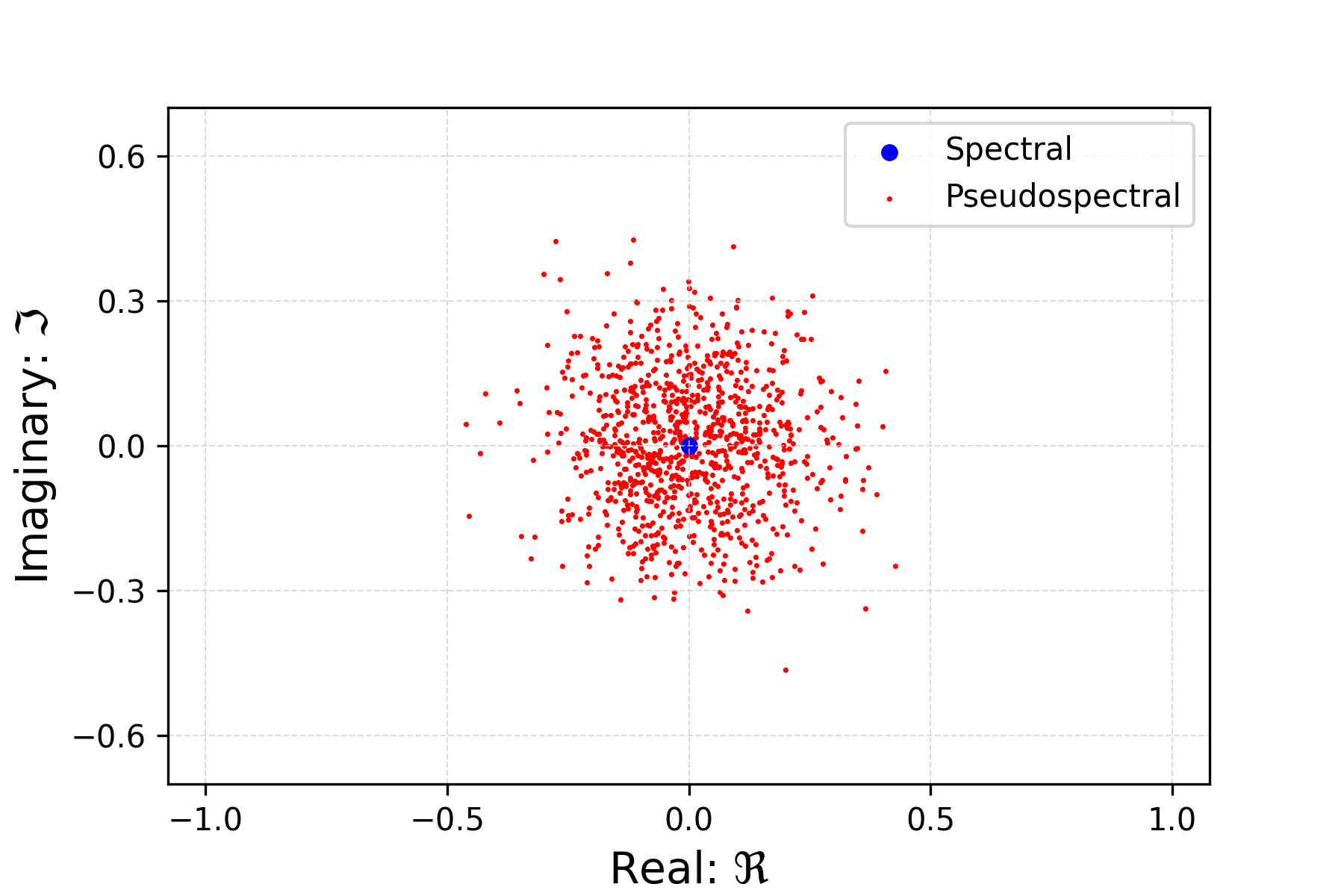}
\caption{$d=100$}
\end{subfigure}
\begin{subfigure}[t]{0.45\linewidth}
\centering
\includegraphics[scale=0.48]{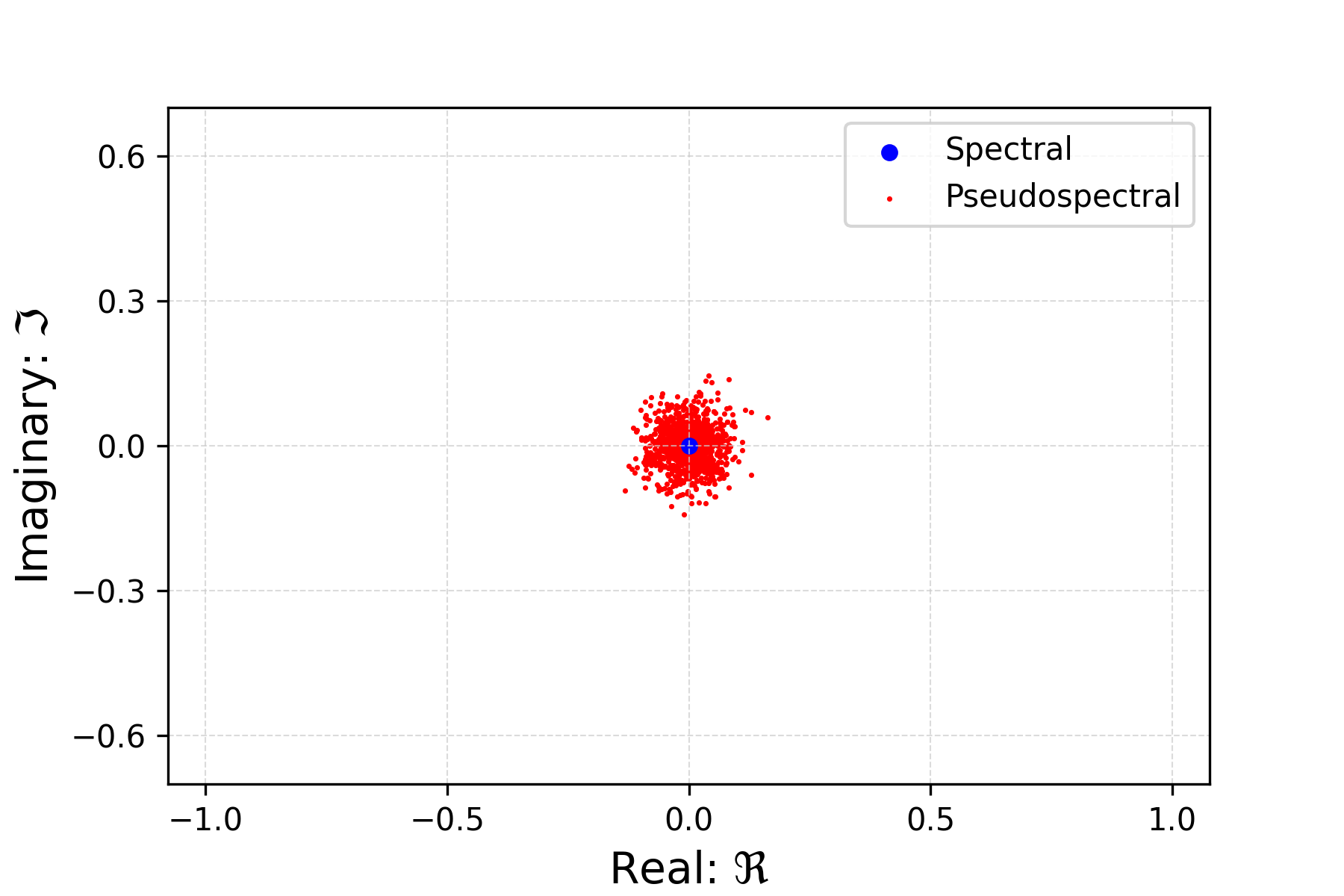}
\caption{$d=1000$}
\end{subfigure}
\caption{Scatter plots illustrating the spectral behavior of rank-1 sampling perturbations. Parameter: $\epsilon = 2$. Number of samples: $N=1000$. } 
\label{fig: zero-matrix}
\end{figure}

We begin our numerical experiments with the simplest case of rank-1 sampling perturbations, where the unperturbed matrix is the zero matrix.  As illustrated in~\Cref{fig: zero-matrix}, the eigenvalues remain tightly concentrated around the origin. This concentration effect becomes more pronounced as the matrix dimension increases. This behavior exemplifies a \textit{regular} perturbation regime in high-dimensional setting: as the dimension tends to infinity, the influence of the perturbation becomes negligible, and the spectrum remains effectively stable. 

This behavior can be explained as follows. According to~\citet[Remark 3.1.1]{vershynin2018high},  the norms of standard normal vectors scale as $\|\pmb{u}\| \sim \sqrt{d}$ and $\|\pmb{v}\|\sim\sqrt{d}$.  By the central limit theorem, the eigenvalues of rank-1 sampling perturbations converge in distribution to a complex normal random variable:
\[
\lambda(E) =  \lambda\left( \frac{ \epsilon \pmb{u}  \pmb{v}^{\dagger}}{\|\pmb{u}\| \|\pmb{v}\|} \right) = \frac{ \epsilon \pmb{v}^{\dagger}\pmb{u}  }{\|\pmb{u}\| \|\pmb{v}\| }\longrightarrow \mathcal{CN}\left( 0, \frac{\epsilon^2}{d} \right), \quad \mathrm{as}\;\;d \rightarrow \infty.
\]
In particular, the effective fluctuation scale shrinks with the standard deviation $\epsilon/\sqrt{d}$. Furthermore, the almost-orthogonality of independent vectors~\citep[Remark 3.2.5]{vershynin2018high} implies the inner product $|\pmb{v}^{\dagger} \pmb{u}| \sim \sqrt{d}$ with high probability, leading to the approximation
\[
\|E\| = \left\| \frac{ \pmb{u}  \pmb{v}^{\dagger}}{\|\pmb{u}\| \|\pmb{v}\|} \right\| = \left| \frac{  \pmb{v}^{\dagger} \pmb{u} }{\|\pmb{u}\| \|\pmb{v}\|} \right| \sim \frac{1}{\sqrt{d}}.
\] 
For a scalar matrix $CI$, every vector in $\mathbb{C}^{d}$ is an eigenvector, and hence the eigenvalues of the perturbed matrix satisfy
\[
\lambda (CI + \epsilon E) = C + \epsilon \lambda (E),
\]
demonstrating that the spectrum remains tightly concentrated around the scalar $C$, with fluctuations on the order of $\epsilon/\sqrt{d}$. This regular concentration behavior extends beyond scalar matrices to all diagonal matrices, or more generally, to normal matrices. A rigorous high-probability analysis of this phenomenon is provided in~\Cref{sec: regular}. 

\subsubsection{Singular pseudospectral concentration}
\label{subsubsec: singular-concentration}

\begin{figure}[htpb!]
\centering
\begin{subfigure}[t]{0.45\linewidth}
\centering
\includegraphics[scale=0.48]{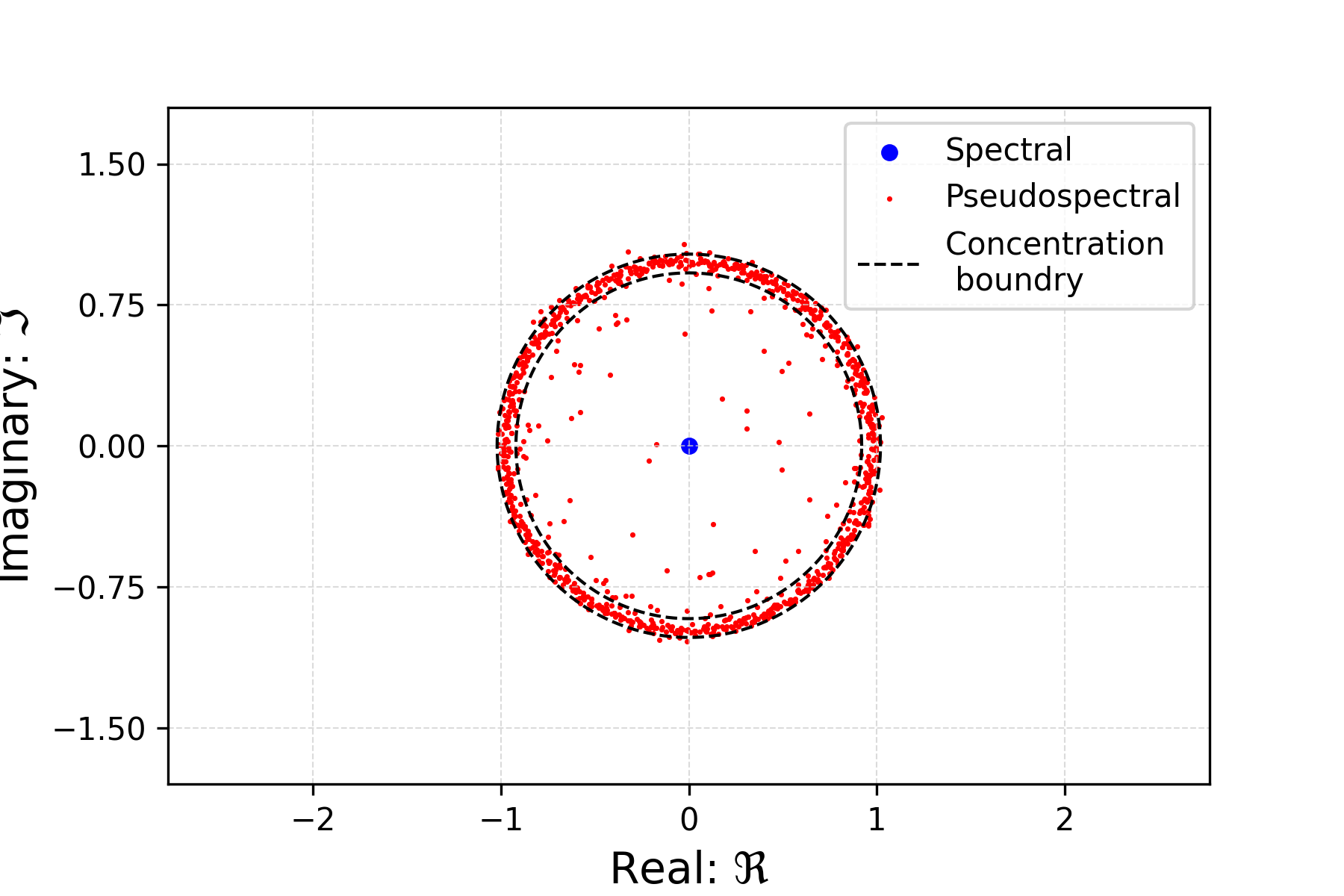}
\caption{$d=100$}
\end{subfigure}
\begin{subfigure}[t]{0.45\linewidth}
\centering
\includegraphics[scale=0.48]{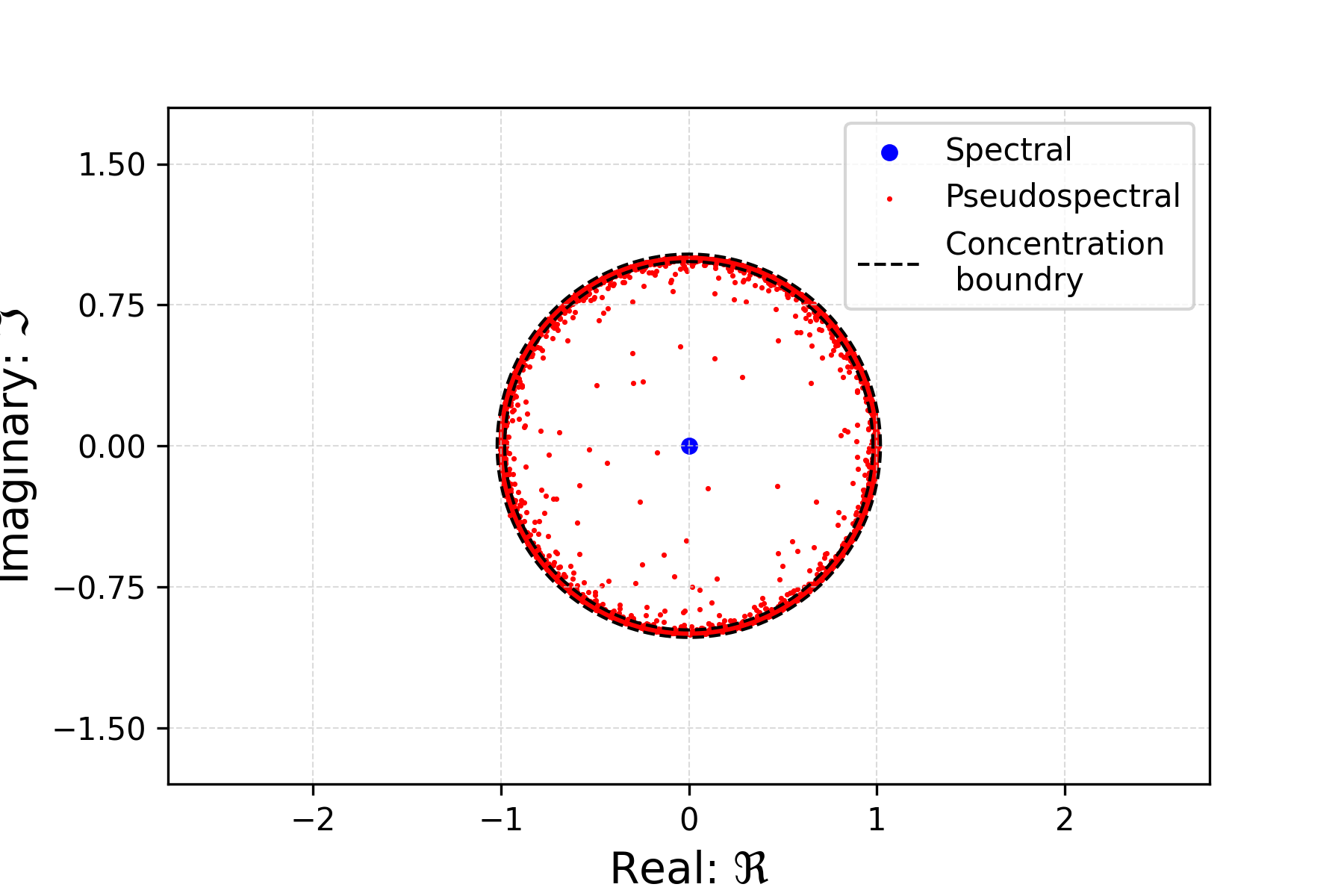}
\caption{$d=1000$}
\end{subfigure}
\caption{Scatter plots illustrating the spectral behavior of the nilpotent Jordan block, as given in~\eqref{eqn: nil-jordan}, under rank-1 sampling perturbations. Parameter: $\epsilon = 2$.  Number of samples: $N=1000$. } 
\label{fig: nilpotent-jordan}
\end{figure}

We now turn to a contrasting example by considering the nilpotent Jordan block, constructed by placing ones on the superdiagonal of the zero matrix:
\begin{equation}
\label{eqn: nil-jordan}
J= 
\begin{pmatrix}
0 & 1 &       &         \\
  & \ddots  & \ddots &   \\
  &          & 0      & 1 \\
  &           &        & 0
\end{pmatrix} \in \mathbb{C}^{d \times d}.
\end{equation}
The corresponding numerical results are illustrated in~\Cref{fig: nilpotent-jordan}.  In contrast to the tightly concentrated eigenvalues near the origin observed in~\Cref{fig: zero-matrix}, the nilpotent Jordan block under rank-1 sampling perturbations exhibits pronounced spectral dispersion: the eigenvalues spread outward, forming an annular region away from the origin. As the matrix dimension increases, the spectral ring becomes narrower and more sharply defined, a striking manifestation of spectral instability.

This stark contrast between the scalar matrix and the nilpotent Jordan block reveals two fundamentally distinct perturbation regimes. In the \textit{regular} regime (\Cref{fig: zero-matrix}), the spectrum remains stable under rank-1 sampling perturbations, with eigenvalues tightly concentrated around the unperturbed values. In the \textit{singular} regime (\Cref{fig: nilpotent-jordan}), however, the same perturbation leads to dramatic spectral displacement and outward dispersion. Although the perturbation norm vanishes in high dimensions, the nilpotent structure significantly amplifies its spectral impact, exhibiting intrinsic spectral sensitivity, a hallmark of \textit{singular perturbation} behavior. 

Previous studies, including~\citet{guionnet2014convergence} and~\citet{basak2020spectrum}, have explored this phenomenon qualitatively, primarily in terms of distributional convergence and weak limits, within the framework of random matrix theory~\citep{tao2012topics}. In contrast, the present work offers a quantitative analysis of this singular perturbation behavior, emphasizing high-probability estimates that more precisely capture the concentration and dispersion behavior induced by random perturbations.

\subsection{A key observation: bottom-left-corner singular perturbation}
\label{subsec: key-observation}

In this section, we highlight a key mechanism underlying the singular perturbation behavior of the nilpotent Jordan block. Specifically, we examine the effect of introducing a single nonzero complex entry $\rho = |\rho|e^{i\theta}$ at the bottom-left corner of the nilpotent Jordan block~\eqref{eqn: nil-jordan}. The resulting perturbed matrix takes the form
\begin{equation}
\label{eqn: nil-jordan-rho}
J(\rho) = 
\begin{pmatrix}
0         & 1               &                   &                \\
           & \ddots       & \ddots        &                \\
           &                  & 0                & 1             \\
\rho     &                  &                   & 0
\end{pmatrix} \in \mathbb{C}^{d \times d}.
\end{equation}
The characteristic polynomial of $J(\rho)$ is given by $\det\left(  \lambda I - J(\rho)  \right) = \lambda^{d} - \rho = 0$.  Hence, the eigenvalues of $J(\rho)$ are the $d$-th roots of $\rho$, explicitly written as: 
\begin{equation}
\label{eqn: nil-jordan-rho-eigenvalue}
\gamma_k =  |\rho|^{\frac{1}{d}} e^{ \frac{ \left(2k\pi  + \theta \right)i }{d} }, \quad k = 0, 1, \ldots, d-1. 
\end{equation}
Thus, all eigenvalues lie uniformly on a circle of radius $|\rho|^{1/d}$ centered at the origin. This is in stark contrast to the unperturbed nilpotent Jordan block, whose eigenvalues are all identically zero. This reveals that even a minimal structural modification, just a single complex entry in the bottom-left corner, induces a complete spectral restructuring.

Importantly, the high-dimensional limit $d \rightarrow \infty$ dramatically amplifies this phenomenon. For instance, when $|\rho| \in \left( (1-\delta)^d, (1+\delta)^d \right)$, then all eigenvalues satisfy $|\gamma_k| \in \left( 1-\delta, 1+\delta \right)$ for $k = 1, \ldots, d$. In the regime where $\delta = \delta_d = C/d$ for some $C > 0$, the annulus $\left( 1-\delta_d, 1+\delta_d \right)$  rapidly collapses onto the unit circle as $d \rightarrow \infty$. Conversely, the range for the perturbation magnitude $|\rho|$, given by $\left( (1-\delta_d)^d, (1+\delta_d)^d \right)$, asymptotically approaches the interval $(e^{-C}, e^{C})$. By choosing $C>0$ sufficiently large, this interval can cover any compact subset of the positive real line $(0, \infty)$, and in the limit $C \rightarrow \infty$, it effectively spans the entire positive real line. This asymptotic behavior reveals a striking high-dimensional spectral sensitivity: even when the perturbation magnitude $|\rho|$ is exponentially small or  large in $d$, the eigenvalues remain tightly concentrated near the unit circle. Thus, the spectral effect of a single bottom-left-corner perturbation becomes delocalized at the origin in large dimensions, producing a nearly uniform distribution of eigenvalues around the unit circle, regardless of the scale of the perturbation.

 

%

Building on this observation, we can generalize the analysis by replacing the fixed entry $\rho$ in the perturbed nilpotent Jordan block~\eqref{eqn: nil-jordan-rho} with a random variable $\xi$ (e.g., $\xi \in \mathcal{CN}(0,1)$). Since the eigenvalues satisfy $|\gamma_k|= |\xi|^{\frac1d}$ for $k = 1, \ldots, d$, the probability of their concentration in the annulus $(1- \delta, 1+ \delta)$ is given by: 
\begin{align}
\mathbf{Pr}\left( \left| | \gamma_k | - 1 \right| < \delta \right) &  = \mathbf{Pr}\left( (1 - \delta)^d <|\xi| < (1 + \delta)^d \right)             \nonumber \\
                                                                                                 & = 1 - \mathbf{Pr}\left( |\xi| \leq (1 - \delta)^d \right) -  \mathbf{Pr}\left( |\xi| \geq (1 + \delta)^d \right). \label{eqn: high-prob-annulus}
\end{align}
For any random variable $\xi$ with bounded density and finite variance, the probability in~\eqref{eqn: high-prob-annulus} approaches unity, with the tail probabilities decaying exponentially in $\delta d$. Setting $\delta = \delta_d = C/d$ for some $C>0$, we obtain the asymptotic estimate
\begin{equation}
\label{eqn: high-prob-annuls-tail}
\lim_{d \rightarrow \infty} \mathbf{Pr}\left( \left| | \gamma_k | - 1 \right| < \frac{C}{d} \right) =  1 - \mathbf{Pr}\left( |\xi| \leq e^{-C} \right) -  \mathbf{Pr}\left( |\xi| \geq e^{C}\right). 
\end{equation}
Under the boundedness assumptions, both tail probabilities $\mathbf{Pr}\left( |\xi| \leq e^{-C} \right)$ and $\mathbf{Pr}\left( |\xi| \geq e^{C}\right)$ decay exponentially. These probabilistic bounds motivate us to propose pseudospectral analysis of the nilpotent Jordan block~\eqref{eqn: nil-jordan-rho} via rank-1 sampling perturbations, which are  rigorously developed in~\Cref{sec: nil-jordan}.




%
%
%
%


We further observe that the bottom-left-corner perturbation is probably the most singular perturbation for the Jordan block~\eqref{eqn: nil-jordan-rho}. Notably, this perturbation remains highly impactful even when the diagonal entries of the Jordan block are identical.  To illustrate this, consider a complex matrix $K \in C^{2d \times 2d}$ with two distinct diagonal values of equal multiplicity $d$. Upon applying a bottom-left-corner perturbation, the perturbed matrix takes the following form: 
\begin{equation}
\label{eqn: two-bottom-left-corner}
K(\rho) = \begin{array}{c}
\begin{pmatrix}
\begin{array}{cccc:cccc}
\gamma_1   & 1            &                      &                      &                      &              &                     &                           \\   
                    & \ddots    & \ddots           &                      &                      &              &                     &                           \\
                    &               & \gamma_1    & 1                   &                      &              &                     &                          \\
                    &               &                      & \gamma_1    & 1                   &              &                     &                          \\  \hdashline            
                    &               &                      &                      & \gamma_2    & 1           &                     &                          \\
                    &               &                      &                      &                      & \ddots   & \ddots           &                          \\
                    &               &                      &                      &                      &              & \gamma_2    &  1                      \\
         \rho     &               &                      &                      &                      &              &                      & \gamma_2
\end{array}
\end{pmatrix}\\
\underbrace{\hspace{3cm}}_{d} \underbrace{\hspace{3cm}}_{d}
\end{array},
\end{equation}
where $\rho \in \mathbb{C}$  represents the perturbation in the bottom-left corner. The characteristic polynomial of $K(\rho)$ takes the form:
\[
\det\left( \lambda I - K(\rho)  \right)   = \left[\lambda^2 - (\gamma_1 + \gamma_2)\gamma + \gamma_1\gamma_2\right]^d - |\rho| e^{\theta i} = 0,
\]
and its eigenvalues are given explicitly by the $d$-th roots of $\rho$, 
\begin{equation}
\label{eqn: two-roots-bottom-left-corner}
\lambda_{k}^{1,2} = \frac{\gamma_1 + \gamma_2 \pm \sqrt{ (\gamma_1 - \gamma_2 )^2 + 4|\rho|^{\frac{1}{d}} e^{ \frac{ \left(2k\pi +  \theta \right)i }{d} }   }   }{2} \quad k = 0, 1, \ldots, d-1. 
\end{equation}
From the expresion~\eqref{eqn: two-roots-bottom-left-corner}, it is evident that when $\gamma_1$ and $\gamma_2$ are well separated, the perturbation exerts only a moderate effect on the eigenvalues.  However, as $\gamma_1$ and $\gamma_2$ approach each other, the influence of the perturbation becomes increasingly significant. The most singular behavior arises in the degenerate case when $\gamma_1 = \gamma_2$, where the spectrum becomes most sensitive to the perturbation. 


\subsection{Overview of contributions and organization}
\label{subsec: overview-contribution}

Based on the observation in~\Cref{subsec: unitary} that random sampling is invariant under unitary similarity, we reduce the analysis of the pseudospectra of complex matrices to the case of upper triangular matrices. Building on the experimental phenomena described in~\Cref{subsec: pseudospectral-phenomenon} and the key mathematical insight of bottom-left-corner singular perturbations from~\Cref{subsec: key-observation}, we develop a concentration theory for the pseudospectra of complex matrices. This leads to the following main contributions:

\begin{itemize}
\item[($\mathbf{I}$)] We first employ moment generating function techniques to derive a tail bound for quadratic forms of standard complex normal vectors. By combining this result with Rouch\'e's theorem from complex analysis, we establish a sub-exponential concentration inequality for the regular pseudospectral concentration of diagonal matrices. Moreover, we show that, with high probability, the separation radius exhibits a dimension-dependent scaling of the form $\delta_d \sim \frac{1}{\sqrt{d}}$, thereby quantifying the rate of regular concentration.
 
\item[($\mathbf{II}$)]  We then generalize the Hanson–Wright concentration inequality and the Carbery–Wright anti-concentration inequality to the complex setting. Exploiting the nilpotent structure of Jordan blocks, and combining this with Rouch\'e's theorem and the Maximum Modulus Principle from complex analysis,  we derive sub-exponential tail bounds for two critical events: escaping beyond the outer radius and falling within the inner radius. This leads to a quantitative theoretical framework for the singular pseudospectral concentration, in which we show that, with high probability, the separation radius again scales with dimension as $\delta_d \sim \frac{1}{\sqrt{d}}$.

\item[($\mathbf{III}$)] Finally, considering upper triangular Toeplitz matrices, which can be expressed via the polynomial symbol of a nilpotent Jordan block, we apply partial fraction decomposition to the associated resolvent, excluding the case of multiple roots. By establishing a root separation bound and using the minimum modulus of the roots to control the convergence of the Laurent series, we extend the theoretical framework for singular pseudospectral concentration to encompass general upper triangular Toeplitz matrices. In this setting as well, we prove that the separation radius scales satisfies $\delta_d \sim \frac{1}{\sqrt{d}}$. 
\end{itemize}

The paper is organized as follows.~\Cref{sec: prelim} introduces the necessary preliminaries from numerical linear algebra, complex analysis, and probabilistic inequalities.~\Cref{sec: regular} develops the theory of regular pseudospectral concentration for normal matrices.~\Cref{sec: nil-jordan}  establishes the theory of singular concentration behavior for nilpotent Jordan blocks, laying out key analytical frameworks.  This analysis is extended in~\Cref{sec: toeplitz}  to general upper-triangular Toeplitz matrices, demonstrating how our theoretical results apply to broader matrix classes. Finally,~\Cref{sec: conclusion} concludes the paper and outlines potential directions for future research.



\section{Preliminaries}
\label{sec: prelim}

This section lays the mathematical foundation for our subsequent analysis by assembling key definitions and theorems across several domains. We begin with tools from pseudospectral theory and matrix analysis,  proceed through relevant results from complex analysis, and conclude with probabilistic inequalities essential for our theoretical developments.

\subsection{Pseudospectra and matrix classes}
\label{subsec: pseudospectra-matrix-classes}

We begin by recalling the definition of pseudospectra, following the second equivalent characterization provided in~\citet[Section 2]{trefethen2020spectra}.
\begin{defn}[Pseudospectra]
\label{defn: pseudospectra}
Let $A \in \mathbb{C}^{d \times d}$ be a matrix, and let $\epsilon>0$ be a small real parameter. The $\epsilon$-pseudospectrum of $A$, denoted $\sigma_{\epsilon}(A)$, is defined as 
\[
\sigma_{\epsilon}(A) = \big\{ \lambda \in \mathbb{C} \big | \lambda \in \sigma(A+\epsilon E)\; \mathrm{for\;some}\; E \in \mathbb{C}^{d \times d}\; \mathrm{with} \; \| E \| = 1 \big\},
\]
where $\sigma(\cdot)$ denotes the spectrum (i.e., the set of eigenvalues). 
\end{defn}

As shown in~\citet[Section 2]{trefethen2020spectra}, the perturbation matrix $E$ may, without loss of generality, be taken to be of rank one.  Let $\pmb{v}$ be a normalized eigenvector of $A + \epsilon E$ corresponding to the eigenvalue $\lambda$, i.e.,
\[
(A + \epsilon E) \pmb{v} = \lambda \pmb{v}.
\]
Suppose further that $E\pmb{v} = s\pmb{u}$, where $\pmb{u}$ is also a normalized vector and $s \in \mathbb{C}$.  Then one can always construct a rank-1 perturbation matrix 
\[
E_1 = s\pmb{u}\pmb{v}^{\top},
\]
which satisfies the same eigenvalue equation:
\[
(A + \epsilon E_1) \pmb{v} = \lambda \pmb{v}.
\] 
We formalize this observation below. 
\begin{theorem}[Rank-1 Characterization of Pseudospectra]
\label{thm: rank1-pseudospectra}
Under the conditions of~\Cref{defn: pseudospectra},  a complex number $\lambda \in \mathbb{C}$ belongs to the $\epsilon$-pseudospectrum of $A$, that is,  $\lambda \in \sigma_{\epsilon}(A)$ if and only if there exists a rank-1 matrix $E \in \mathbb{C}^{d \times d}$ with $\| E \| = 1$ such that $\lambda \in \sigma(A + \epsilon E)$.
\end{theorem}

To facilitate the presentation of classical results, we next introduce several classes of matrices. Let $\mathscr{C}(d)$ denote the space of all complex $d \times d$ matrices. The subset $\mathscr{N}(d) \subset \mathscr{C}(d)$ represents the set of normal matrices, and its complement $\mathscr{C}(d) \setminus \mathscr{N}(d)$ consists of nonnormal matrices. Among important subclasses, we denote by $\mathscr{U}(d)$ the set of unitary matrices, $\mathscr{R}(d)$ the set of upper triangular matrices, and $\mathscr{D}(d)$ the set of diagonal matrices. We also consider the class of Toeplitz matrices $\mathscr{T}(d)$, and in particular their upper-triangular counterparts, referred to as $\mathcal{R}$-Toeplitz matrices and denoted by $\mathscr{T}(R, d)$.  

According to the Schur decomposition, any complex matrix is unitarily similar to an upper-triangular matrix, and any normal matrix is unitarily similar to a diagonal matrix. These classical results are summarized below for completeness.
\begin{lemma}[Schur Decomposition, Theorem 7.1.3 in~\citet{golub2013matrix}]
\label{lem: schur}
For any matrix $A \in \mathscr{C}(d)$, there exists a unitary matrix $U \in \mathscr{U}(d)$ such that 
\[
U^{\dagger}AU = T,
\]
where $T \in \mathscr{T}(d)$ is an upper-triangular matrix.
\end{lemma}

\begin{lemma}[Corollary 7.1.4 in~\citet{golub2013matrix}]
\label{lem: schur-diag}
For any normal matrix $N \in \mathscr{N}(d)$, there exists a unitary matrix $U \in \mathscr{U}(d)$ such that 
\[
U^{\dagger}NU = D
\]
where $D =\mathrm{diag}(\gamma_1, \ldots, \gamma_d) \in \mathscr{D}(d)$ is a diagonal matrix.
\end{lemma}
As discussed in~\Cref{subsec: unitary}, the statistical properties of both full-rank and rank-1 sampling perturbations are invariant under unitary similarity. In particular, based on~\eqref{eqn: unitary-similar-full-rank} ---~\eqref{eqn: perturb-unitary}, we may analyze the pseudospectral behavior of a matrix under such random perturbations via its Schur decomposition.
\begin{theorem}(Pseudospectral Invariance under Schur Decomposition)
\label{thm: triangle}
Let $A \in \mathscr{C}(d)$, and let $T \in \mathscr{R}(d)$ be the upper-triangular matrix obtained via the Schur decomposition, as described in~\Cref{lem: schur}. Then, the pseudospectral behavior of $A$ under both full-rank and rank-1 sampling perturbations is identical to that of $T$. Moreover, if $A \in \mathscr{N}(d)$, its pseudospectral behavior coincides with that of the diagonal matrix $D=\mathrm{diag}(\gamma_1, \ldots, \gamma_d) \in \mathscr{D}(d)$, as described in~\Cref{lem: schur-diag}. 
\end{theorem}

We also require a classical result concerning rank-1 updates of determinants. 
\begin{lemma}[Matrix Determinant Lemma, Section 2.1.4 in~\citet{golub2013matrix} ]
\label{lem: matrix-determinant}
Let $\pmb{u}, \pmb{v} \in \mathbb{C}^d$ be two column vectors, and let  $A \in \mathscr{C}(d)$ be an invertible matrix. Then, 
\begin{equation}
\label{thm: matrix-det}
\det(A + \pmb{u} \pmb{v}^{\dagger}) = \det(A) \left( 1 + \pmb{v}^{\dagger}A^{-1} \pmb{u} \right).
\end{equation}
\end{lemma}

\subsection{Tools from complex analysis}
\label{subsec: tool-complex-analysis}

To study the roots of the characteristic polynomial $\det(\lambda I - A - \epsilon E)$ and relate them to those of the unperturbed matrix $\det(\lambda I - A)$, we invoke classical results from complex analysis. 
\begin{lemma}[Rouch\'e's Theorem, Theorem 4.3 of Chapter 3 in~\citet{stein2010complex}]
\label{lem: rouche}
Suppose that $f$ and $g$ are holomorphic in an open set containing a circle $\Gamma$ and its interior. If
\[
\left| f(z) \right| > \left| g(z) \right| \quad \mathrm{for\; all} \; z \in \Gamma,
\] 
then $f$ and $f + g$ have the same number of zeros inside the circle $\Gamma$, counted with multiplicities.
\end{lemma}

\begin{lemma}[Maximum Modulus Principle, Theorem 4.5 of Chapter 3 in~\citet{stein2010complex}]
\label{lem: maximum-modulus-principle}
If $f$ is a non-constant holomorphic function in a connected open region $\Omega$, then $|f(z)|$ cannot attain a maximum at any point in the interior of $\Omega$; that is, the maximum of $|f(z)|$ must occur on the boundary of $\Omega$. 
\end{lemma}


\subsection{Probabilistic inequalities}
\label{subsec: prob-inq-prelim}

To support our probabilistic analysis, we now introduce several  fundamental inequalities concerning real standard normal vectors. 

\begin{lemma}[Hanson-Wright Concentration Inequality,~\citet{hanson1971bound}]
\label{lem: hanson-wright}
Let $\pmb{X} \sim \mathcal{N}(\pmb{0}, I)$ be a standard normal vector in $\mathbb{R}^d$, and let $A \in \mathbb{R}^{d \times d}$ be a real matrix. Then, there exists a universal constant $C>0$ such that for any $t > 0$,
\begin{equation}
\label{eqn: hanson-wright}
\mathbf{Pr}\left( \big| \pmb{X}^{\top}A\pmb{X} - E[\pmb{X}^{\top}A\pmb{X}] \big| \geq t \right) \leq 2 \exp \left( -C \min \left( \frac{t^2}{\|A\|_F^2}, \frac{t}{\|A\|} \right)\right),
\end{equation}
where $\|A\|_{F}^{2} = \sum\limits_{i=1}^{d} \sum\limits_{j=1}^{d} a_{ij}^2$ is the Frobenius norm.
\end{lemma}

\begin{lemma}[Carbery–Wright Anti-Concentration Inequality, Theorem 8 in~\citet{carbery2001distributional} ]
\label{lem: carbery-wright}
Let $\pmb{X}  \sim \mathcal{N}(\pmb{0}, I)$ be a standard normal vector in $\mathbb{R}^d$, and let $f: \mathbb{R}^d \mapsto \mathbb{R}$ be a real-valued polynomial  of degree at most $m \leq d$. Then, there exists a universal constant $C>0$ such that for any $t>0$,
\begin{equation}
\label{eqn: carbery-wright}
\mathbf{Pr}\left( \left| f(\pmb{X}) \right| \leq t  \sqrt{  \mathrm{Var}[ f(\pmb{X}) ]  }  \right) \leq Cm t^{\frac1m}.
\end{equation}
\end{lemma}

To avoid reliance on distributional norms for real-valued normal vectors, we conclude with a concentration inequality for the norm of a complex normal vector.
\begin{lemma}[Concentration of Norm]
\label{lem: norm-concentration}
Let $\pmb{\xi} \sim \mathcal{CN}( \pmb{0},I)$  be a complex normal  vector in $\mathbb{C}^d$. Then, there exist a universal constant $C>0$ such that for any $t > 0$,
\begin{equation}
\label{eqn: norm-concentration}
\mathbf{Pr} \big( \left|  \| \pmb{\xi}\|^2 - d  \right| \leq td \big)  \leq  Ce^{- \sqrt{d} t}.
\end{equation}
\end{lemma}
The complete proof is provided in~\Cref{sec: app-1}.



%
%
%

\section{Regular pseudospectral concentration for normal matrices}
\label{sec: regular}

In this section, we investigate the pseudospectral properties of normal matrices under rank-1 sampling perturbations. By~\Cref {thm: triangle}, the analysis can be reduced, without loss of generality, to the case of diagonal matrices: 
\begin{equation}
\label{eqn: diag-matrix}
D = \mathrm{diag}(\gamma_1, \ldots, \gamma_d) = \mathrm{diag}(\gamma_{1,r} + i\gamma_{1,i}, \ldots, \gamma_{d,r} + i\gamma_{d,i})  \in \mathscr{D}(d).
\end{equation}
This simplification remains valid even when the perturbations originate from rank-1 sampling mechanisms. To lay the groundwork for analyzing regular pseudospectral concentration in such diagonal settings, we first derive a concentration inequality for quadratic forms involving complex normal vectors.

\subsection{Tail bound for quadratic forms of complex normal vectors}
\label{subsec: tail-bound}

Let $\pmb{u}, \pmb{v} \in \mathbb{C}^d$ be two complex vectors, where the real and imaginary parts of each component are given as
\[
u_{k} = u_{k,r} + iu_{k,i}, \qquad v_{k} = v_{k,r} + iv_{k,i}.
\]
Expanding the product $\gamma_k \overline{v_k} u_k$, we obtain:
\begin{align*} 
\gamma_k \overline{v_k} u_k = & \left[ \gamma_{k,r} \left( v_{k,r}u_{k,r}+ v_{k,i}u_{k,i} \right) - \gamma_{k,i}\left( v_{k,r} u_{k, i} -  v_{k,i}  u_{k,r}\right) \right]  \\
                                                  & + i\left[ \gamma_{k,r}\left( v_{k,r} u_{k,i} -  v_{k,i}  u_{k,r} \right) + \gamma_{k,i}\left( v_{k,r}u_{k,r} + v_{k,i}u_{k,i} \right) \right].
\end{align*}
Summing over $k=1,\ldots, d$, we decompose the quadratic form $\pmb{v}^{\dagger}D\pmb{u} $ into its real and imaginary parts as
\begin{subequations}
\begin{align}
&  \Re \big( \pmb{v}^{\dagger}D\pmb{u}  \big) =  \sum_{k=1}^{d}   \left[ \gamma_{k,r} \left( v_{k,r}u_{k,r}+ v_{k,i}u_{k,i} \right) - \gamma_{k,i}\left( v_{k,r} u_{k, i} -  v_{k,i}  u_{k,r}\right) \right]  \label{eqn: real-coll}\\
&  \Im \big(\pmb{v}^{\dagger}D\pmb{u}  \big) = \sum_{k=1}^{d} \left[ \gamma_{k,r}\left( v_{k,r} u_{k,i} -  v_{k,i}  u_{k,r} \right) + \gamma_{k,i}\left( v_{k,r}u_{k,r} + v_{k,i}u_{k,i} \right) \right]            \label{eqn: imag-coll}
\end{align}
\end{subequations}
 Assume now that $\pmb{u}, \pmb{v} \sim  \mathcal{CN}(\pmb{0},I)$ are independent standard complex normal vectors. Then, each of the real and imaginary components is an i.i.d. real normal random variable (i.e.,  $u_{k,r}, u_{k,i}, v_{k,r}, v_{k,i} \sim  \mathcal{N}(0,1/2)$). From the symmetry in~\eqref{eqn: real-coll} and~\eqref{eqn: imag-coll}, we observe that the real part $ \Re \left(  \pmb{v}^{\dagger}D\pmb{u}  \right)$ and the imaginary part $\Im \left(  \pmb{v}^{\dagger}D\pmb{u} \right)$ are i.i.d. real random variables. In particular, we have:
\[
\mathbf{Pr}\left( \big| \Re \big( \pmb{v}^{\dagger}D\pmb{u} \big) \big| \geq 2td  \right) = \mathbf{Pr}\left( \big| \Im \big( \pmb{v}^{\dagger}D\pmb{u} \big) \big| \geq 2td  \right). 
\]
Applying the union bound, we obtain:
\begin{align}
\mathbf{Pr}\left( \big| \pmb{v}^{\dagger}D\pmb{u} \big| \geq 2\sqrt{2}td  \right)  & \leq \mathbf{Pr}\left(  \big| \Re \big(  \pmb{v}^{\dagger}D\pmb{u} \big) \big| \geq 2td  \right) + \mathbf{Pr}\left( \big| \Im \big(\pmb{v}^{\dagger}D\pmb{u} \big) \big| \geq 2td  \right) \nonumber \\
                                                                                                                              & = 2\mathbf{Pr}\left(  \big| \Re \big(  \pmb{v}^{\dagger}D\pmb{u} \big) \big| \geq 2td  \right). \label{eqn: prob-complex-real}
\end{align}
We now proceed to derive a tail bound for the real part.
\begin{lemma}
\label{lem: regular-con1}
Let $\pmb{u}, \pmb{v} \sim \mathcal{CN}(\pmb{0}, I)$ be two independent standard complex normal vectors, and let $\delta > 0$ be a small positive constant. Suppose the coefficients $\gamma_k \in \mathbb{C}$ for $k=1,\ldots, d$ satisfy the uniform bound 
\[
\max_{1 \leq k \leq d} |\gamma_k| = \frac{1}{\delta}.
\]
Then, there exists a universal constant $C > 0$ such that for any $t > 0$, the following tail bound holds:
\begin{equation}
\label{eqn: regular-con1}
\mathbf{Pr}\left( \big|  \Re \big( \pmb{v}^{\dagger}D\pmb{u} \big) \big| \geq 2td  \right) \leq C e^{-  t \delta \sqrt{2d}}. 
\end{equation}
\end{lemma}
\begin{proof}[Proof of~\Cref{lem: regular-con1}]
For any $\lambda \in (0, \delta)$, we apply Markov’s inequality to the exponential moment associated with the real part $ \pmb{v}^{\dagger}D\pmb{u}$, as expressed in~\eqref{eqn: real-coll}.  This yields the following upper bound on the tail probability:
\begin{equation}
\label{eqn: regular-tail-1}
\mathbf{Pr}\left( \Re \big(  \pmb{v}^{\dagger}D\pmb{u} \big)\geq 2td  \right)  \leq e^{-2\lambda td} \mathbb{E}\left[ \prod_{k=1}^{d} e^{ 2 \lambda \left[ \gamma_{k,r} \left( v_{k,r}u_{k,r}+ v_{k,i}u_{k,i} \right) - \gamma_{k,i}\left( v_{k,r} u_{k, i} -  v_{k,i}  u_{k,r}\right) \right] } \right].
\end{equation}
Since the components $u_{k,r}, u_{k,i}, v_{k,r}, v_{k,i} \sim  \mathcal{N}(0,1/2)$, for $k=1, \ldots, d$, are i.i.d.~real normal variables for $k=1,\ldots, d$, the expectation factorizes a product of independent terms:
\begin{align}
& \mathbb{E}\left[ \prod_{i=1}^{d} e^{ 2\lambda \left[ \gamma_{k,r} \left( v_{k,r}u_{k,r}+ v_{k,i}u_{k,i} \right) - \gamma_{k,i}\left( v_{k,r} u_{k, i} -  v_{k,i}  u_{k,r}\right) \right] } \right] \nonumber  \\
& \phantom{\mathrel{abscdefghijk}} =   \prod_{i=1}^{d} \mathbb{E}\left[  e^{ 2\lambda  \gamma_{k,r} \left( v_{k,r}u_{k,r}+ v_{k,i}u_{k,i} \right) } \right]  \cdot \mathbb{E} \left[ e^{ 2\lambda \gamma_{k,i}\left(  v_{k,i}  u_{k,r}  - v_{k,r} u_{k, i} \right)} \right] \nonumber \\
& \phantom{\mathrel{abscdefghijk}} = \prod_{i=1}^{d} \left(\mathbb{E}\left[  e^{ 2\lambda  \gamma_{k,r} v_{k,r}u_{k,r} } \right] \right)^2 \cdot \left( \mathbb{E} \left[ e^{ 2\lambda \gamma_{k,i}v_{k,r} u_{k, i} } \right] \right)^2. \label{eqn: regular-con2}
\end{align}
Each term in the expectation~\eqref{eqn: regular-con2} involves the moment generating function of the product of two independent normal variables. For such random variables $u,v \sim \mathcal{N}(0,1/2)$, and any $s \in (0, 1/\delta)$,  a standard computation yields: 
\begin{equation}
\label{eqn: regular-con3}
\mathbb{E}\left[  e^{ 2\lambda  s vu} \right]  = \frac{1}{\pi} \int_{-\infty}^{\infty} \int_{-\infty}^{\infty} e^{2 \lambda s xy - x^2 - y^2} dxdy  = \frac{1}{\sqrt{1 - \lambda^2 s^2 }}.   
\end{equation}
Substituting~\eqref{eqn: regular-con3} into the expectation~\eqref{eqn: regular-con2}, we have: 
\begin{equation}
\label{eqn: regular-con4}
 \mathbb{E}\left[ \prod_{i=1}^{d} e^{ 2\lambda \left[ \gamma_{k,r} \left( v_{k,r}u_{k,r}+ v_{k,i}u_{k,i} \right) - \gamma_{k,i}\left( v_{k,r} u_{k, i} -  v_{k,i}  u_{k,r}\right) \right] } \right] = \prod_{k=1}^d \frac{1}{ \left(1 - \lambda^2 \gamma_{k,r}^2\right) \left(1 - \lambda^2 \gamma_{k,i}^2\right)}.
\end{equation}
Using the identity $|\gamma_k|^2 = \gamma_{k,r}^2 + \gamma_{k,i}^2$, we observe:
\begin{equation}
\label{eqn: regular-con5}
\left( 1 - \lambda^2 \gamma_{k,r}^2 \right)\left( 1 - \lambda^2 \gamma_{k,i}^2 \right) = 1 - \lambda^2 |\gamma_k|^2 + \lambda^4 \gamma_{k,r}^2  \gamma_{k,i}^2  \geq 1 - \frac{\lambda^2}{\delta^2}, 
\end{equation}
where the inequality follows from the assumption $|\gamma_k| \leq 1/\delta$ for $k=1,\ldots, d$. Substituting~\eqref{eqn: regular-con3} and~\eqref{eqn: regular-con5} into the tail bound~\eqref{eqn: regular-tail-1}, we obtain:
\[
\mathbf{Pr}\left( \Re \big(  \pmb{v}^{\dagger}D\pmb{u}  \big)\geq 2td  \right)  \leq \exp\left( -2 \lambda  td -  \log\left(1 - \frac{\lambda^2}{ \delta^2}\right) d \right).
\]
Applying the inequality $\log(1 - x) \geq -x - x^2 $ for any $x \in (0,1/2)$, we further obtain: 
\begin{equation}
\label{eqn: regular-con6}
\mathbf{Pr}\left( \Re \big( \pmb{v}^{\dagger}D\pmb{u} \big)\geq 2td  \right)  \leq \exp\left(- \left( 2 \lambda t - \frac{\lambda^2}{\delta^2} - \frac{\lambda^4}{\delta^4} \right)d \right),
\end{equation}
which holds for any $\lambda \in (0, \delta/\sqrt{2})$. Finally, by choosing $\lambda = \delta / \sqrt{2d}$ in the tail bound~\eqref{eqn: regular-con6} and setting $C = e^{\frac34}$, we complete the proof. 
%
\end{proof}

According to the tail inequality established in~\eqref{eqn: prob-complex-real}, a similar concentration bound holds for the modulus $|\pmb{v}^{\dagger}D\pmb{u}|$, up to a factor of $2$. The formal statement is given below.
\begin{lemma}
\label{lem: regular-con2}
Under the same assumptions as in~\Cref{lem: regular-con1}, there exists a universal constant $C > 0$ such that for any $t > 0$, the following tail bound holds:
\begin{equation}
\label{eqn: regular-con21}
\mathbf{Pr}\left( \big|   \pmb{v}^{\dagger}D\pmb{u}  \big| \geq  td  \right) \leq C e^{- \frac{t \delta \sqrt{d}}{2} }. 
\end{equation}
\end{lemma}
\subsection{Regular concentration behavior}
\label{subsec: regular-concentration}

\begin{figure}[htpb!]
\centering
\begin{subfigure}[t]{0.45\linewidth}
\centering
\includegraphics[scale=0.48]{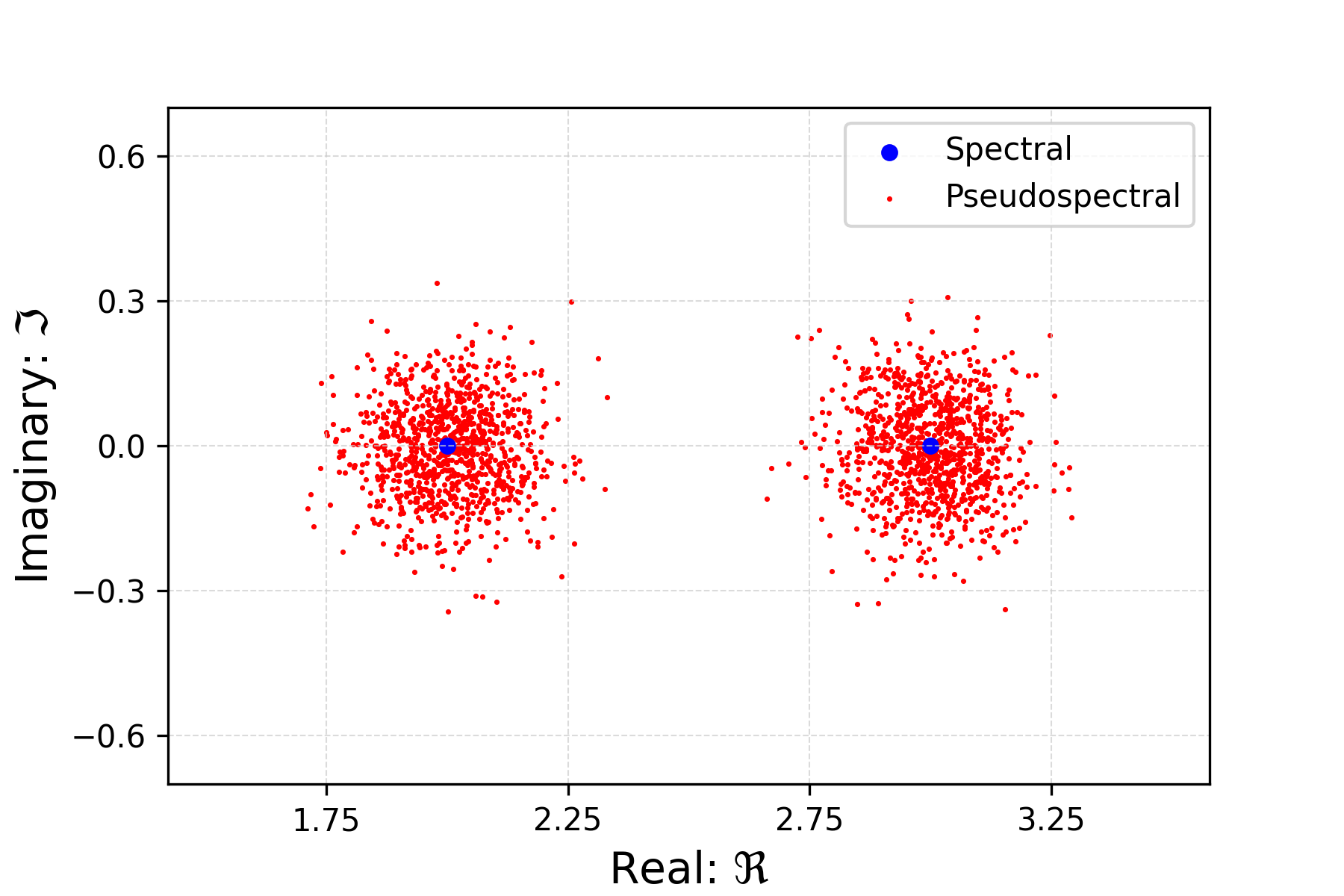}
\caption{$d=100$}
\end{subfigure}
\begin{subfigure}[t]{0.45\linewidth}
\centering
\includegraphics[scale=0.48]{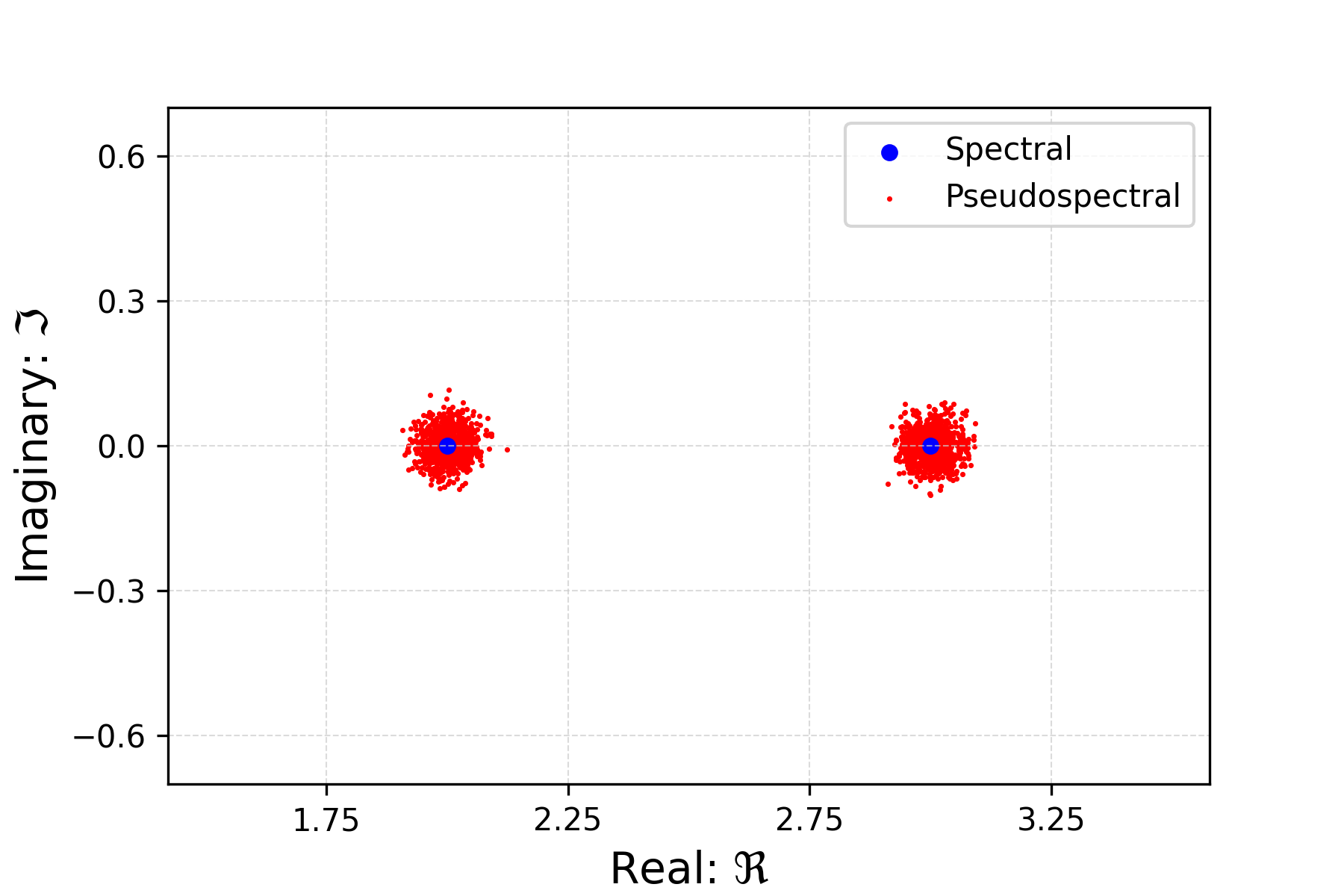}
\caption{$d=1000$}
\end{subfigure}
\caption{Scatter plots illustrating the spectral behavior of a diagonal matrix under rank-1 sampling perturbations. The underlying diagonal matrix contains two distinct eigenvalues, $2$ and $3$, each repeated with multiplicity $d/2 $. Parameter: $\epsilon = 2$.  Number of samples: $N=1000$.} 
\label{fig: diagonal-bi}
\end{figure}

As described in~\Cref{subsubsec: regular-concentration},  the pseudospectra of the scalar matrix under rank-1 sampling perturbations form a circular disk centered at the scalar value. Notably, as the matrix dimension increases from $d=100$ to $d=1000$, the radius of the disk contracts by approximately a factor of $\sqrt{10}$. To  explore this behavior more generally, we consider a diagonal matrix 
\[
D =\mathrm{diag}(\gamma_1, \ldots, \gamma_d)  \in \mathscr{D}(d),
\] 
with distinct entries, i.e., $\gamma_k\neq \gamma_{\ell}$ for $k \neq \ell$, where $k,\ell =1,\ldots, d$. The corresponding numerical experiment is presented in~\Cref{fig: diagonal-bi}. In this setup, the diagonal matrix $D \in \mathscr{D}(d)$ has two distinct values, each appearing with multiplicity $d/2$. The resulting pseudospectral plots exhibit a structural pattern similar to that in~\Cref{fig: zero-matrix}, with clusters of eigenvalues forming around the original spectrum.

Building upon the empirical observations from~\Cref{fig: diagonal-bi}, we now proceed to rigorously characterize the eigenvalue distribution under rank-1 sampling perturbations. As a preliminary step, we introduce some necessary definitions. Let the separation radius be defined as
\[
\delta_0 := \frac{1}{2} \min_{k\neq \ell }\big |\gamma_k - \gamma_{\ell} \big |, \qquad \mathrm{for}\quad k,\ell = 1, 2, \ldots, d,
\]
which guarantees that for any $\delta \leq \delta_0$, the open neighborhoods 
\[
O(\gamma_k; \delta) := \left\{ \lambda \in \mathbb{C}: |\lambda - \gamma_k| < \delta \right\}
\]  
are mutually disjoint.  The union of these disjoint neighborhoods defines the set
\begin{equation}
\label{eqn: disjoint-disk}
S(\delta) := \bigcup\limits_{k = 1}^{d} O(\gamma_k; \delta) = \bigcup\limits_{k = 1}^{d}\left\{ \lambda \in \mathbb{C} :  |\lambda - \gamma_k| < \delta \right\},
\end{equation}
which consists of $d$ non-intersecting disks centered at the eigenvalues of the diagonal matrix $D =\mathrm{diag}(\gamma_1, \ldots, \gamma_d)$. We now state the main result concerning regular spectral concentration under rank-1 sampling perturbations.
\begin{theorem}[Regular Spectral Concentration]
\label{thm: regular-concentration}
Let $D =  \mathrm{diag}(\gamma_1, \ldots, \gamma_d) \in \mathscr{D}(d)$ be a diagonal matrix, and let $\pmb{u}, \pmb{v} \sim \mathcal{CN}(\pmb{0},I)$ be two independent standard complex normal vectors. For any $\epsilon > 0$, there exists universal constants $d_0 > 0$ and $C>0$ such that, for any dimension $d \geq d_0$, the following concentration inequality holds:
\begin{equation}
\label{eqn: regular-concentration}
\mathbf{Pr}\left(  \sigma\left( D +  \frac{  \epsilon  \pmb{u}  \pmb{v}^{\dagger}}{\|\pmb{u}\| \|\pmb{v}\|} \right) \subset S(\delta_d) \right)\geq 1 - \epsilon,
\end{equation}
where the set $S(\delta_d)$ is defined in~\eqref{eqn: disjoint-disk}, and the separation radius is given by 
\[
\delta_d = \frac{4C\epsilon}{\sqrt{d}}.
\]
\end{theorem}

\begin{proof}[Proof of~\Cref{thm: regular-concentration}]
For any $\delta > 0$, observe that when $\lambda \in \partial S(\delta)$, the matrix $\lambda I - D$ is invertible. We now proceed to analyze the characteristic polynomial of the perturbed matrix using~\Cref{lem: matrix-determinant} (Matrix Determinant Lemma). For a normalized rank-1 perturbation, we obtain the exact factorization:
\begin{equation}
\label{eqn: diagonal-1}
\det\left( \lambda I - D - \frac{  \epsilon  \pmb{u}  \pmb{v}^{\dagger}}{\|\pmb{u}\| \|\pmb{v}\|} \right) = \det\left( \lambda I - D  \right) \cdot \left( 1 - \frac{\epsilon \pmb{v}^{\dagger} (\lambda I -D)^{-1} \pmb{u} }{ \|\pmb{v}\| \|\pmb{u}\| } \right).
\end{equation}
Next, by applying~\Cref{lem: rouche} (Rouch\'{e}'s Theorem), we deduce that if the following condition holds for any $\lambda \in  \partial S(\delta)$:
\begin{equation}
\label{eqn: diagonal-2}
\frac{\epsilon \left| \pmb{v}^{\dagger} (\lambda I -D)^{-1} \pmb{u} \right|  }{ \|\pmb{v}\| \|\pmb{u}\| }  < 1,
\end{equation}
then the number of roots of the perturbed characteristic polynomial inside each disk is preserved. In other words, the eigenvalues of the perturbed matrix remain within the union $S(\delta)$ of disjoint neighborhoods centered at the original eigenvalues.

From the factorization in~\eqref{eqn: diagonal-1} and the condition in~\eqref{eqn: diagonal-2}, we derive the following probabilistic bound as
\begin{align}
\mathbf{Pr}\bigg( \sigma\left( D +  \frac{  \epsilon  \pmb{u}  \pmb{v}^{\dagger}}{\|\pmb{u}\| \|\pmb{v}\|} \right) \subset S(\delta) \bigg)   & \geq \mathbf{Pr}\left( \sup_{\lambda \in \partial S(\delta)} \frac{ \left| \epsilon \pmb{v}^{\dagger} (\lambda I -D)^{-1} \pmb{u} \right|  }{ \|\pmb{v}\| \|\pmb{u}\| }  < 1 \right)   \nonumber \\
                                                                     & =  1 -  \mathbf{Pr}\left(  \sup_{\lambda \in \partial S(\delta)}  \bigg| \sum_{k=1} \frac{\overline{v}_k u_k}{\lambda - \gamma_k} \bigg| \geq \frac{1}{ \epsilon}  \|\pmb{u}\| \|\pmb{v}\| \right).          \label{eqn: diagonal-3}                
\end{align}
To further bound this probability, we apply the union bound:
\begin{multline}
\mathbf{Pr}\left(  \sup_{\lambda \in \partial S(\delta)}  \bigg| \sum_{k=1} \frac{\overline{v}_k u_k}{\lambda - \gamma_k} \bigg|  \geq \frac{1}{ \epsilon}  \|\pmb{u}\|  \|\pmb{v}\| \right) \leq  \mathbf{Pr}\left(\|\pmb{u}\|^2 <(1-t)d \right) + \mathbf{Pr}\left(\|\pmb{v}\|^2 <(1-t)d \right)   \\ 
 + \mathbf{Pr}\left(  \sup_{\lambda \in \partial S(\delta)}  \bigg| \sum_{k=1} \frac{\overline{v}_k u_k}{\lambda - \gamma_k} \bigg|  \geq \frac{1}{ \epsilon}  \|\pmb{u}\|  \|\pmb{v}\|,\; \|\pmb{u}\|^2 \geq (1-t)d, \; \|\pmb{v}\|^2 \geq (1-t)d \right),      \label{eqn: diagonal-4} 
\end{multline}
for any $t \in (0,1)$. The constants appearing in the tail bounds from~\eqref{eqn: norm-concentration} and~\eqref{eqn: regular-con21} are universal, so we denote $C_1 > 0$ as their maximum. Choose a constant $C>0$ such that $3C_1 e^{-C} = \epsilon$. By applying~\Cref{lem: norm-concentration} with $t= \frac{C}{\sqrt{d}}$, we obtain the following bound as
\begin{equation}
\label{eqn: diagonal-5} 
\mathbf{Pr}\left(\|\pmb{u}\|^2 <\left(1-   \frac{C}{\sqrt{d}} \right)d \right) + \mathbf{Pr}\left(\|\pmb{v}\|^2 <\left(1-   \frac{C}{\sqrt{d}} \right) d \right) \leq 2C_1 e^{-C} = \frac{2\epsilon}{3}.
\end{equation}
For the remaining term in the union bound~\eqref{eqn: diagonal-4}, we employ conditional probability to refine the estimate as follows:
\begin{align}
& \mathbf{Pr}\left(  \sup_{\lambda \in \partial S(\delta)}  \bigg| \sum_{k=1} \frac{\overline{v}_k u_k}{\lambda - \gamma_k} \bigg|  \geq \frac{1}{ \epsilon}  \|\pmb{u}\|  \|\pmb{v}\|,\; \|\pmb{u}\|^2 \geq \left(1-   \frac{C}{\sqrt{d}} \right) d, \; \|\pmb{v}\|^2 \geq \left(1-   \frac{C}{\sqrt{d}} \right) d \right)  \nonumber \\
 &  \phantom{\mathrel{+}} \leq \mathbf{Pr}\left(   \sup_{\lambda \in \partial S(\delta)}  \bigg| \sum_{k=1} \frac{\overline{v}_k u_k}{\lambda - \gamma_k} \bigg|  \geq \frac{1}{ \epsilon}  \|\pmb{u}\|  \|\pmb{v}\| \bigg |\; \|\pmb{u}\|^2 \geq \left(1-   \frac{C}{\sqrt{d}} \right) d, \; \|\pmb{v}\|^2 \geq \left(1-   \frac{C}{\sqrt{d}} \right)d \right)  \nonumber \\
 & \phantom{\mathrel{+}} \leq \mathbf{Pr}\left(  \sup_{\lambda \in \partial S(\delta)}  \bigg| \sum_{k=1} \frac{\overline{v}_k u_k}{\lambda - \gamma_k} \bigg|  \geq \frac{1}{ \epsilon} \left(1-   \frac{C}{\sqrt{d}} \right) d \right). \label{eqn: diagonal-6} 
\end{align}
By the compactness of $\partial S(\delta)$,  there exists a point $\lambda_0 \in \partial S(\delta)$ that achieves the supremum:
\[
 \sup_{\lambda \in \partial S(\delta)}  \bigg| \sum_{k=1} \frac{\overline{v}_k u_k}{\lambda - \gamma_k} \bigg| = \bigg| \sum_{k=1} \frac{\overline{v}_k u_k}{\lambda_0 - \gamma_k} \bigg|.
\]
Since there exists $d_0>0$ such that for any $d \geq d_0$, the inequality $1 - \frac{C}{\sqrt{d}} \geq \frac{1}{2}$ holds. Thus, we apply~\Cref{lem: regular-con2} to derive the exponential tail bound as
\begin{align}
&\mathbf{Pr}\left(  \sup_{\lambda \in \partial S(\delta)}  \bigg| \sum_{k=1} \frac{\overline{v}_k u_k}{\lambda - \gamma_k} \bigg|  \geq  \frac{1}{ \epsilon} \left(1 - \frac{C}{\sqrt{d}}\right)d \right)  \leq \mathbf{Pr}\left( \bigg| \sum_{k=1} \frac{\overline{v}_k u_k}{\lambda_0 - \gamma_k} \bigg|  \geq \frac{d}{ 2\epsilon } \right)  \leq C_1 e^{- \frac{\delta\sqrt{d} }{4\epsilon} }.  \label{eqn: diagonal-7} 
\end{align}
By summarizing inequalities~\eqref{eqn: diagonal-3} through~\eqref{eqn: diagonal-7}, and choosing the separation radius $\delta = \delta_d = \frac{4C\epsilon}{\sqrt{d}}$, we conclude the proof. 

\end{proof}

\section{Singular pseudospectral concentration for nilpotent Jordan block}
\label{sec: nil-jordan}

As illustrated in~\Cref{fig: nilpotent-jordan},  the nilpotent Jordan block defined in~\eqref{eqn: nil-jordan}, when subjected to rank-1 sampling perturbations, exhibits a striking phenomenon known as singular pseudospectral concentration. As discussed in~\Cref{subsubsec: singular-concentration}, the eigenvalues of the perturbed matrix spread outward, forming an annular region rather than concentrating near the origin. Remarkably, the width of this annulus diminishes as the matrix dimension increases. In this section, we provide a rigorous mathematical analysis of this singular behavior.

\subsection{Two key probabilistic inequalities for complex normal vectors}
\label{subsec: complex-probabilstic-ineq}
 
Before analyzing the pseudospectral behavior, we extend the two key probabilistic inequalities,~\eqref{eqn: hanson-wright} and~\eqref{eqn: carbery-wright}, from the real to the complex setting.  These results serve as preparatory tools for the forthcoming analysis.  Complete proofs are provided in~\Cref{sec: technical-detail-complex-ineq}. We begin with a complex variant of the Hanson–Wright concentration inequality (cf.~\Cref{lem: hanson-wright}), which provides tail bounds for quadratic forms involving independent complex normal vectors.

\begin{lemma}[Complex Hanson-Wright Concentration Inequality]
\label{lem: hanson-wright-complex}
Let $\pmb{u}, \pmb{v} \sim \mathcal{CN}(\pmb{0}, I)$ be two independent standard complex normal vectors in $\mathbb{C}^d$. Then there exists a universal constant $C>0$ such that for any $t > 0$ and any integer $k \in \{0, \ldots, d-1\}$, the following inequality holds:
\begin{equation}
\label{eqn: complex-hanson-wright}
\mathbf{Pr}\left( \big| \pmb{v}^{\dagger}J_0^k \pmb{u}  \big| \geq t \right) \leq 4 \exp \left( -C \min \left( \frac{t^2}{d - k}, t \right)\right),
\end{equation} 
\end{lemma}
 
As a counterpart to~\Cref{lem: carbery-wright}, we now introduce a complex analog of the Carbery–Wright anti-concentration inequality, which quantifies the probability that quadratic forms involving complex normal vectors are close to zero.
\begin{lemma}[Complex Carbery–Wright Anti-Concentration Inequality]
\label{lem: carbery-wright-complex}
Let $\pmb{u}, \pmb{v} \sim \mathcal{CN}(\pmb{0}, I)$ be two independent standard complex normal vectors in $\mathbb{C}^d$, and let $A \in \mathscr{C}(d)$ be a complex matrix. Then there exists a universal constant $C>0$ such that for any $t>0$, the following inequality holds:
\begin{equation}
\label{eqn: complex-carbery-wright}
\mathbf{Pr}\left( \big| \pmb{v}^{\dagger} A \pmb{u} \big| \leq t  \|A\|_{F}    \right) \leq C \sqrt{t}.
\end{equation}
\end{lemma}
 
\subsection{Singular concentration behavior}
\label{subsec: singular-concentration}

Since the Jordan block $J$, defined in~\eqref{eqn: nil-jordan}, is nilpotent, that is, $J^d = 0$, the resolvent $(I - a J)^{-1}$ admits a closed-form expression for any $a \in \mathbb{C}$. This follows from the fact that the corresponding Neumann series terminates after $d$ terms due to the nilpotency of $J$. We formalize this observation in the following lemma.
\begin{lemma}
\label{lem: inverse-jordan}
Let $J \in \mathscr{C}(d)$ be the nilpotent Jordan block given in~\eqref{eqn: nil-jordan},  then $J^d = 0$. Furthermore, for any $a \in \mathbb{C}$, the resolvent $(I - a J)^{-1}$ is given explicitly by 
\begin{equation}
\label{eqn: inverse-jorda}
(I - a J)^{-1} = \sum_{k=0}^{d-1} a^{k} J^{k}.
\end{equation}
\end{lemma}

Before analyzing the spectrum of the nilpotent Jordan block $J$ under rank-1 sampling perturbations,  we first evaluate the probability that the unperturbed eigenvalue $\lambda = 0$ remains an eigenvalue after perturbation. Since the first column of $J$ contains only zeros and elementary column operations preserve the rank, the invertibility of the perturbed matrix
\[
J + \frac{  \epsilon  \pmb{u}  \pmb{v}^{\dagger}}{\|\pmb{u}\| \|\pmb{v}\| }
\]
is equivalent  to that of a matrix where only the first column is modified, that is,
\[
J +  \left( \frac{  \epsilon \pmb{u}\overline{v}_1}{\|\pmb{u}\| \|\pmb{v}\| }, \pmb{0}, \ldots, \pmb{0} \right). 
\]
In this representation, the determinant is nonzero if and only if the lower-left corner entry, specifically, the product $u_d \overline{v}_1$, is nonzero. Since $u_d$ and $\overline{v}_1$ are independent standard complex normal random variables (i.e., $u_d, \overline{v}_1 \sim \mathcal{CN}(0,1)$), their product is nonzero almost surely. Hence, we conclude that
\begin{equation}
\label{eqn: lambda-neq-a}
\mathbf{Pr}(\lambda = 0) = \mathbf{Pr}\left( \mathrm{rank}\left(J + \frac{  \epsilon  \pmb{u}  \pmb{v}^{\dagger}}{\|\pmb{u}\| \|\pmb{v}\|}  \right) = d-1 \right) =0.  
\end{equation}
This guarantees that $\lambda = 0$ is almost surely not an eigenvalue of the perturbed matrix. Thus,  we can apply~\Cref{lem: matrix-determinant} (Matrix Determinant Lemma) to express the characteristic polynomial as
\begin{align*}
\det\left( \lambda I - J - \frac{\epsilon \pmb{u} \pmb{v}^{\dagger}}{\| \pmb{u} \| \| \pmb{v} \|} \right) &  = \det \left( \lambda I - J \right) \cdot  \left(1  -  \frac{\epsilon \pmb{v}^{\dagger} (\lambda I - J)^{-1} \pmb{u}  }{\| \pmb{v} \| \| \pmb{u} \| }  \right)                                                               \nonumber  \\                                                                                                             
                                                                                                                                                         & =  \lambda^{d}  \left(1  -  \frac{\epsilon \pmb{v}^{\dagger} \left( I - \frac{J}{\lambda} \right)^{-1} \pmb{u}  }{\lambda \| \pmb{v} \| \| \pmb{u} \| }  \right).  
\end{align*}                                                                                                                                                      
Next, applying~\Cref{lem: inverse-jordan},  which provides an explicit expansion of the resolvent, we further derive the following closed-form expression for the characteristic polynomial as 
\begin{equation}
\label{eqn: jordan-polynomial-2} 
\det\left( \lambda I - J - \frac{\epsilon \pmb{u} \pmb{v}^{\dagger}}{\| \pmb{u} \| \| \pmb{v} \|} \right) = \lambda^{d}  \left( 1 - \frac{\epsilon}{\| \pmb{v} \| \| \pmb{u} \| } \sum_{k=0}^{d-1}\frac{ \pmb{v}^{\dagger}J^k \pmb{u} }{\lambda^{k+1}} \right).
\end{equation}

 We now present a rigorous statement for the phenomenon of singular spectral concentration for the nilpotent Jordan block, as described in~\eqref{eqn: nil-jordan}. 
\begin{theorem}[Singular Spectral Concentration, Nilpotent Jordan Block]
\label{thm: nilpotent-jordan}
Let $J \in \mathscr{C}(d)$ be the nilpotent Jordan block defined in~\eqref{eqn: nil-jordan}, and let $\pmb{u}, \pmb{v} \sim \mathcal{CN}(0,1)$ two independent standard complex normal vectors. For any $\epsilon > 0$, there exists universal constants $d_0 > 0$ and $C>0$ such that, for any dimension $d \geq d_0$, the following concentration inequality holds:
\begin{equation}
\label{eqn: nilpotent-jordan-concentration}
\mathbf{Pr}\left( \sigma\left( J +  \frac{  \epsilon  \pmb{u}  \pmb{v}^{\dagger}}{\|\pmb{u}\| \|\pmb{v}\|} \right)   \subseteq \big\{ \lambda \in \mathbb{C} : \big |  | \lambda |  - 1\big | < \delta_d  \big\} \right)  \geq 1 - \epsilon,
\end{equation}
where the concentration radius satisfies 
\[
\delta_d = 2\sqrt{\frac{C\epsilon }{d}}.
\] 
\end{theorem}
\begin{proof}[Proof of~\Cref{thm: nilpotent-jordan}]
From~\eqref{eqn: lambda-neq-a}, we know that almost surely, $\lambda = a_0$ is not an eigenvalue of the perturbed matrix. To analyze the probability that all eigenvalues lie within the annular region, we consider the complement of this event. This complement consists of two disjoint events: eigenvalues either escape beyond the outer radius or fall within the inner radius. Accordingly, we write:
\begin{align}
& \mathbf{Pr}\left( \sigma\left( J +  \frac{  \epsilon  \pmb{u}  \pmb{v}^{\dagger}}{\|\pmb{u}\| \|\pmb{v}\|} \right) \subseteq \big\{ \lambda \in \mathbb{C}:  \big | | \lambda| - 1 \big| < \delta \big\} \right)  \nonumber \\
 &  = 1 - \underbrace{\mathbf{Pr}\left( \exists  \lambda \in \mathbb{C} \bigg |   \det\left( \lambda I - J -  \frac{  \epsilon  \pmb{u}  \pmb{v}^{\dagger}}{\|\pmb{u}\| \|\pmb{v}\|} \right) = 0 \quad \mathrm{and }\quad  | \lambda | \geq 1 + \delta \right)}_{:=\mathbf{I}}  \nonumber \\
 &  \mathrel{\phantom{= 1}} - \underbrace{\mathbf{Pr}\left( \exists  \lambda \in \mathbb{C} \bigg |   \det\left( \lambda I - J -  \frac{  \epsilon  \pmb{u}  \pmb{v}^{\dagger}}{\|\pmb{u}\| \|\pmb{v}\|} \right) = 0 \quad \mathrm{and }\quad  0< | \lambda | \leq 1 - \delta  \right)}_{:=\mathbf{II}}. \label{eqn: linear-symbol-thm-1}
\end{align}

\begin{itemize}

\item[(1)] \textbf{Escaping beyond the outer radius.} We begin by analyzing the probability that an eigenvalue escapes beyond the outer boundary of the annular region. Along this boundary, the spectral parameter satisfies
\[
\lambda  = (1+\delta)e^{i\theta},\quad \mathrm{where}\;\;\theta \in [0, 2\pi).
\]
To quantify the likelihood that the rank-1 sampling perturbation drives an eigenvalue beyond this contour, we invoke~\Cref{lem: rouche} (Rouch\'{e}'s Theorem). In particular, we  aim to estimate the escape probability
\begin{equation}
\label{eqn: linear-symbol-outer-1}
\mathbf{I}  \leq \mathbf{Pr} \left( \epsilon \left| \sum_{k=0}^{d-1}\pmb{v}^{\dagger} \Omega_k \pmb{u} \right| \geq  \| \pmb{v} \| \| \pmb{u} \|  \right),
\end{equation}
where the perturbation terms are defined as
\[
\Omega_k := \frac{J^k }{(1+\delta)^{k+1} e^{i(k+1)\theta}}, \quad k=0, 1, \ldots, d-1. 
\]
Given that $\pmb{u}, \pmb{v} \sim \mathcal{CN}(\pmb{0},I)$ are two independent standard complex normal vectors in $\mathbb{C}^d$, we decompose this event via the union bound. For any $t>0$, we have:
\begin{align}
\mathbf{I} & \leq  \mathbf{Pr} \left( \| \pmb{u} \|^2 < (1-t)d \right) + \mathbf{Pr} \left( \| \pmb{v} \|^2 < (1-t)d \right)  \nonumber \\ 
                 & \mathrel{\phantom{\leq}}+  \mathbf{Pr} \left( \epsilon \left|\sum_{k=0}^{d-1}\pmb{v}^{\dagger} \Omega_k \pmb{u} \right| \geq \| \pmb{u} \| \| \pmb{v} \|, \| \pmb{u} \|^2 \geq (1-t)d,  \| \pmb{v} \|^2 \geq (1-t)d \right)  \nonumber \\
                 & \leq  \mathbf{Pr} \left( \| \pmb{u} \|^2 < (1-t)d \right) + \mathbf{Pr} \left( \| \pmb{v} \|^2 < (1-t)d \right)  +  \mathbf{Pr} \left( \epsilon \left|\sum_{k=0}^{d-1}\pmb{v}^{\dagger} \Omega_k \pmb{u} \right|\geq  (1-t)d \right).          \label{eqn: linear-symbol-outer-2}
\end{align}
To estimate the first two terms, we apply~\Cref{lem: norm-concentration}. Let $C_1 > 0$ be the universal constant in the tail inequality~\eqref{eqn: norm-concentration}, and choose $C_2 > 0$ such that $6C_1e^{C_2} = \epsilon$. Then, setting $t= \frac{C_2}{\sqrt{d}}$, we obtain: 
\begin{equation}
\label{eqn: linear-symbol-outer-3}
 \mathbf{Pr} \left( \| \pmb{u} \|^2 < \left(1- \frac{C_2}{\sqrt{d}} \right)d \right) +  \mathbf{Pr} \left( \| \pmb{v} \|^2 < \left(1- \frac{C_2}{\sqrt{d}} \right)d \right) \leq \frac{\epsilon}{3}. 
\end{equation}

Next, we turn to the third term in the union bound~\eqref{eqn: linear-symbol-outer-2}. By applying the triangle inequality and bounding the modulus, we obtain the following estimate:      
\begin{equation}
\mathbf{Pr} \left( \epsilon \left| \sum_{k=0}^{d-1}\pmb{v}^{\dagger} \Omega_k \pmb{u} \right| \geq  \left(1- \frac{C_2}{\sqrt{d}} \right) d \right) \leq \mathbf{Pr} \left( \epsilon \sum_{k=0}^{d-1}  \left|\frac{  \pmb{v}^{\dagger}J^k \pmb{u} }{(1+\delta)^{k+1} }\right| \geq  \left(1- \frac{C_2}{\sqrt{d}} \right) d \right).  \label{eqn: linear-symbol-outer-4}
\end{equation} 
To proceed, we estimate the relevant geometric sum: 
\begin{equation}
\label{eqn: linear-symbol-outer-5}
\sum_{k=0}^{d-1} \frac{d-k}{(1+\delta)^{k+1}}
= \frac{1}{1+\delta} \cdot \frac{1 - (d+1)(1+\delta)^{-d} + d(1+\delta)^{-d-1}}{\left(1 - \frac{1}{1+\delta}\right)^2} \leq \frac{1+\delta}{\delta^2}.
\end{equation}
Since there exists a universal constant $d_{0,1} >0$ such that for any $d\geq d_{0,1}$, it holds that $1-C_2/\sqrt{d} \geq 1/2$,  we may relax the right-hand side by replacing the threshold with $d/2$. Substituting this bound~\eqref{eqn: linear-symbol-outer-5} into~\eqref{eqn: linear-symbol-outer-4}, we obtain: 
\begin{align}
\mathbf{Pr} \left( \epsilon \left| \sum_{k=0}^{d-1} \Omega_k \right| \geq   \left(1- \frac{C_2}{\sqrt{d}} \right)d \right) & \leq \mathbf{Pr} \left( \epsilon \left| \sum_{k=0}^{d-1} \Omega_k \right| \geq  \frac{d}{2} \right)  \nonumber \\
                                                                                                                                     & \leq  \sum_{k=0}^{d-1}\mathbf{Pr} \left(  \big|  \pmb{v}^{\dagger}J^k \pmb{u} \big| \geq \frac{ (d-k) \delta^2 d}{2\epsilon (1+\delta)}  \right).  \label{eqn: linear-symbol-outer-6}
\end{align} 

Assuming the concentration radius satisfies
\[
\delta_d = 2\sqrt{\frac{ C\epsilon }{ d}},
\]
for some constant $C>0$ to be determined later. Then, there exists a threshold $d_{0,2} >0$ such that for any $d\geq d_{0,2}$, the inequality $1+\delta \leq 2$ holds. Under this assumption, we can simplify the bound in~\eqref{eqn: linear-symbol-outer-6} as follows: 
 \begin{equation}
 \label{eqn: linear-symbol-outer-7}
\mathbf{Pr} \left( \epsilon \left| \sum_{k=0}^{d-1} \Omega_k \right| \geq   \left(1- \frac{C_2}{\sqrt{d}} \right) d \right)  \leq  \sum_{k=0}^{d-1}\mathbf{Pr} \left(  \big|  \pmb{v}^{\dagger}J^k \pmb{u} \big| \geq  (d-k)C  \right). 
\end{equation}
To estimate the right-hand side, we invoke~\Cref{lem: hanson-wright-complex} (Complex Hanson-Wright Concentration Inequality), which provides a tail bound for complex quadratic forms. Let $C_3>0$ denote the universal constant in the bound given by~\eqref{eqn: complex-hanson-wright}. By applying~\Cref{lem: hanson-wright-complex} into each term in the sum~\eqref{eqn: linear-symbol-outer-7}, we obtain: 
 \begin{equation}
 \label{eqn: linear-symbol-outer-8}
\mathbf{Pr} \left( \epsilon \left| \sum_{k=0}^{d-1} \Omega_k \right| \geq  (1-t)d \right)  \leq  4\sum_{k=0}^{d-1}\exp\left( -(d-k)C_3 \min\left( C^2, C \right) \right)
\end{equation}
To ensure the total escape probability is sufficiently small, we choose $C>0$ such that 
\[
\exp\left( -C_3 \min\left( C^2, C \right) \right) = \min\left\{ \frac{1}{2}, \frac{\epsilon}{48} \right\}. 
\]
This choice ensures that the total contribution of the tail probabilities satisfies $\mathbf{I} \leq \epsilon/2$. 


\item[(2)] \textbf{Falling within the inner radius.} We now consider the probability that an eigenvalue falls within the inner boundary of the annular region.  Along this boundary, the spectral parameter satisfies
\[
\lambda  = (1-\delta) e^{i\theta},\quad \mathrm{where}\;\; \theta \in [0, 2\pi).
\]
To estimate the likelihood that the rank-1 sampling perturbation pushes an eigenvalue within this contour, we invoke~\Cref{lem: maximum-modulus-principle} (Maximum Modulus Principle). Specifically, we bound the probability of such inward deviation as follows:
\begin{equation}
\label{eqn: linear-symbol-inner-1}
\mathbf{II}  \leq \mathbf{Pr} \left( \epsilon \left| \sum_{k=0}^{d-1} \pmb{v}^{\dagger} \omega_k \pmb{u} \right| \leq  \| \pmb{v} \| \| \pmb{u} \|  \right),
\end{equation}
where the perturbation terms are given by
\[
\omega_k := \frac{J^k }{(1-\delta)^{k+1} e^{i(k+1)\theta}}, \quad k=0, 1, \ldots, d-1. 
\]
Applying the union bound, we can decompose the probability in~\eqref{eqn: linear-symbol-inner-1} as follows: 
\begin{equation}
\label{eqn: linear-symbol-inner-2}
\mathbf{II} \leq   \mathbf{Pr} \left( \| \pmb{u} \|^2 > (1+t)d \right) + \mathbf{Pr} \left( \| \pmb{v} \|^2 > (1+ t)d \right) + \mathbf{Pr} \left( \epsilon \left| \sum_{k=0}^{d-1} \pmb{v}^{\dagger} \omega_k \pmb{u}  \right| \leq  (1+t)d \right). 
\end{equation}
The first two terms are bounded via~\Cref{lem: norm-concentration}.  Since we set $t= \frac{C_2}{\sqrt{d}}$, where $C_2 > 0$ is chosen so that $6C_1e^{C_2} = \epsilon$, we obtain the bound as
\begin{equation}
\label{eqn: linear-symbol-inner-3}
\mathbf{Pr} \left( \| \pmb{u} \|^2 > \left(1+\frac{C_2}{\sqrt{d}} \right)d \right) + \mathbf{Pr} \left( \| \pmb{v} \|^2 > \left(1+\frac{C_2}{\sqrt{d}} \right)d \right) \leq \frac{\epsilon}{3}. 
\end{equation}
For the final term in~\eqref{eqn: linear-symbol-inner-2}, we invoke~\Cref{lem: carbery-wright-complex} (Complex Carbery–Wright Anti-Concentration Inequality), which require us to estimate the Frobenius norm as 
\begin{equation}
\label{eqn: linear-symbol-inner-4}
\left\| \sum_{k=0}^{d-1} \omega_k \right\|_{F}^{2} = \sum_{k=0}^{d-1} \frac{d-k}{(1-\delta)^{2(k+1)}} \geq \frac{1}{(1-\delta)^{2d}}.
\end{equation}
As shown earlier, the concentration radius satisfies
\[
\delta_d = 2\sqrt{\frac{C \epsilon }{ d}},
\]
where $C> 0$ is a constant defined previously. By applying the bound for the Frobenius norm from~\eqref{eqn: linear-symbol-inner-4}, we know that there exists a constant $d_{0,3}>0$ such that for all $d \geq d_{0,3}$, 
\begin{equation}
\label{eqn: linear-symbol-inner-5}
\left\| \sum_{k=0}^{d-1} \omega_k \right\|_{F} \geq \frac{1}{(1-\delta_d)^{d}} \geq  \frac{\epsilon }{ 36C_4^2}\left( 1 + \frac{C_2}{\sqrt{d}} \right)d,
\end{equation}
where $C_4>0$ is  the universal constant from the tail inequality~\eqref{eqn: complex-carbery-wright}. This ensures that  the final term in~\eqref{eqn: linear-symbol-inner-2} is at most $\epsilon/6$, and therefore the falling probability satisfies $\mathbf{II} \leq \epsilon/2$.

\end{itemize}

Finally, by choosing $d_0 = \max\{d_{0,1}, d_{0,2}, d_{0,3}\}$, we ensure that all previous bounds hold uniformly for any $d \geq d_0$, thereby completing the proof. 
\end{proof}

\section{Application to $\mathcal{R}$-Toeplitz matrices}
\label{sec: toeplitz}

In this section, we study the class of $\mathcal{R}$-Toeplitz matrices $\mathscr{T}(R,d)$, and analyze their pseudospectral behavior under rank-1 sampling perturbations. Each $\mathcal{R}$-Toeplitz matrix $T \in \mathscr{T}(R,d)$ is uniquely determined by an associated polynomial symbol of degree at most $n$ (where $1 \leq n \leq d-1$), given by:
\begin{equation}
\label{eqn: symbol}
 p(z) = \sum_{k=0}^{n} a_k z^k, 
\end{equation}
where the constant term $a_0$ is referred to as the trivial coefficient, and the remaining coefficients $a_k$ for $k=1,\ldots,n$ are considered nontrivial. The corresponding $\mathcal{R}$-Toeplitz matrix realization takes the canonical form:
\begin{equation}
\label{eqn: ut-toeplitz}
T =  T(a_0, a_1, \ldots, a_{n}) 
= \begin{pmatrix} 
    a_0   & a_1         & \cdots        & a_{n}               &               &                   \\
             & \ddots     & \ddots        & \ddots              & \ddots    &                   \\
             &                & \ddots        & \ddots              & \ddots    & a_{n}         \\
             &                &                   & \ddots              & \ddots    & \vdots        \\
             &                &                   &                         & \ddots    & a_1           \\
             &                &                   &                         &               & a_0 
   \end{pmatrix}        \in \mathscr{T}(R,d),
\end{equation}
where each $k$-th superdiagonal is filled with the constant value $a_k$ from the polynomial symbol~\eqref{eqn: symbol}. A fundamental special case is the nilpotent Jordan block defined in~\eqref{eqn: nil-jordan}, which corresponds to the linear symbol with $p(z)= z$, where $a_0 =0$ and $a_1 =1$. This observation highlights an important structural connection between $\mathcal{R}$-Toeplitz matrices and polynomial functions of nilpotent operators. As a result, the theoretical framework developed in~\Cref{sec: nil-jordan} can be naturally extended to this broader class of structured matrices.

\subsection{Binomial-symbol $\mathcal{R}$-Toeplitz matrices}
\label{subsec: binomial}

\begin{figure}[htpb!]
\centering
\begin{subfigure}[t]{0.45\linewidth}
\centering
\includegraphics[scale=0.48]{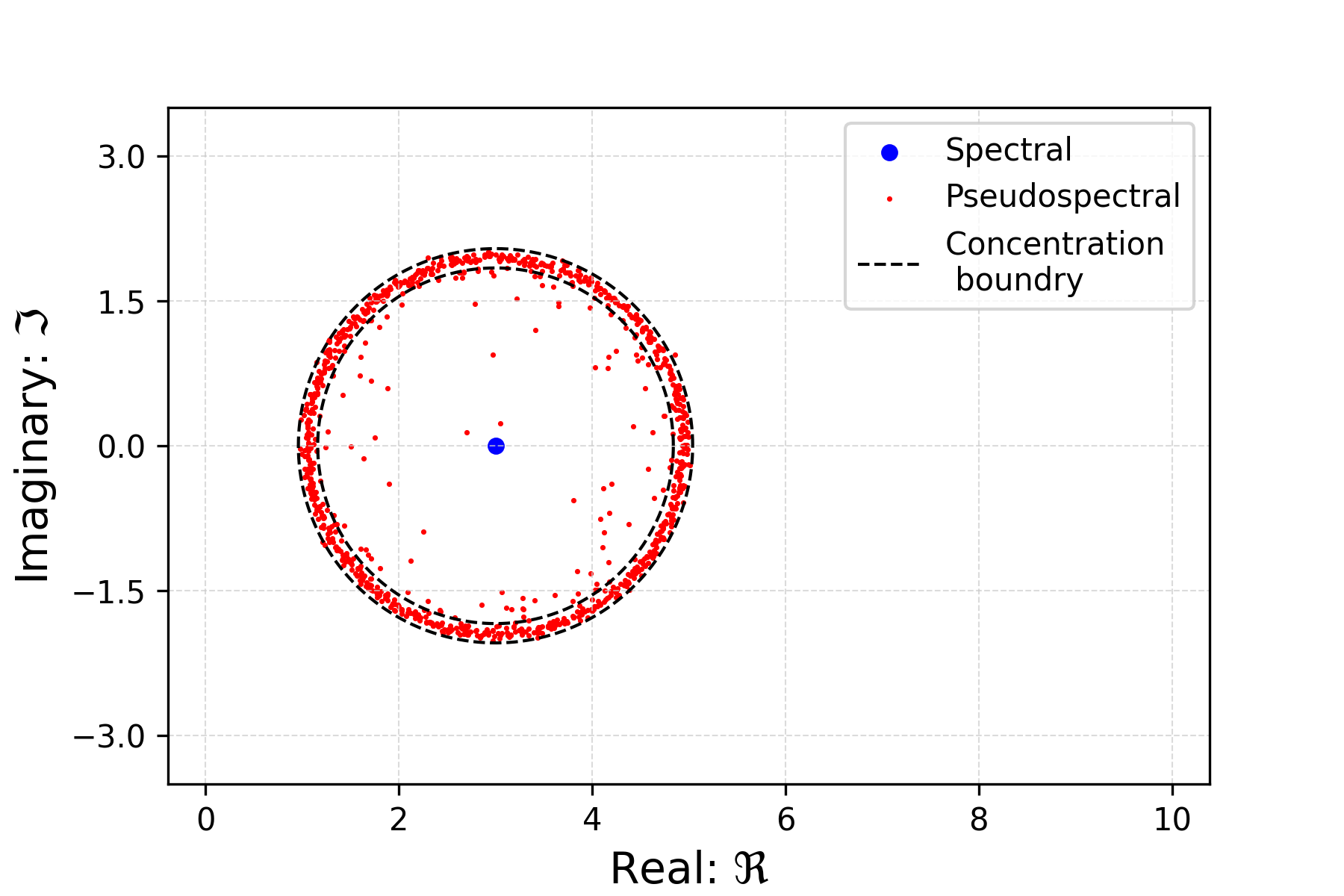}
\caption{$d=100$}
\end{subfigure}
\begin{subfigure}[t]{0.45\linewidth}
\centering
\includegraphics[scale=0.48]{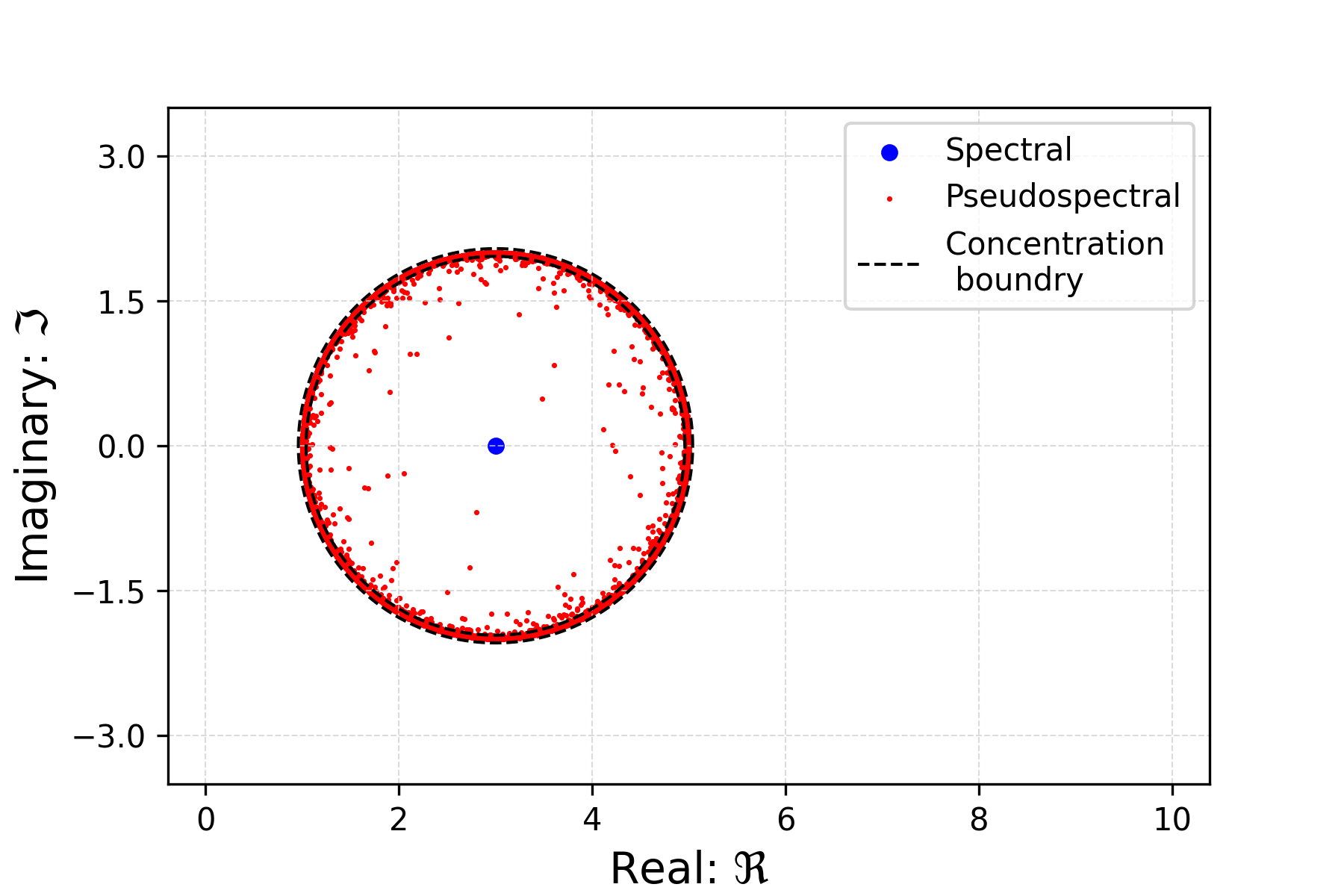}
\caption{$d=1000$}
\end{subfigure}
\caption{Scatter plots illustrating the spectral behavior of the $\mathcal{R}$-Toeplitz matrix defined by a linear symbol $p(z)=3 + 2z$ under rank-1 sampling perturbations. Parameter: $\epsilon = 2$.  Number of samples: $N=1000$.} 
\label{fig: linear-symbol}
\end{figure}

Our analysis begins with $\mathcal{R}$-Toeplitz matrices generated by binomial symbols, where only one coefficient beyond the constant term is nontrivial. Specifically, the symbol takes the form
\begin{equation}
\label{eqn: binomial}
p(z) = a_mz^m + a_0,
\end{equation}
where $m \in \{1,\ldots, n\}$. A simple yet illustrative case arises when the symbol is linear, i.e., $p(z) = a_0 + a_1 z$. In this case, the corresponding $\mathcal{R}$-Toeplitz matrix is given by
\begin{equation}
\label{eqn: linear-symbol-teoplitz}
T = T(a_0,a_1) = a_0I + a_1J, 
\end{equation}
where $I$ is the identity matrix and $J$ is the standard nilpotent Jordan block. As defined in~\eqref{eqn: nil-jordan}, the nilpotent Jordan block corresponds to the special case where $a_0=0$ and $a_1 = 1$. Building upon the singular concentration behavior observed for the nilpotent Jordan block in~\Cref{fig: nilpotent-jordan}, we now examine the spectral behavior of an $\mathcal{R}$-Toeplitz matrix generated by the linear symbol $p(z)=3+2z$ under rank-$1$ sampling perturbations, as shown in~\Cref{fig: linear-symbol}.  Compared to~\Cref{fig: nilpotent-jordan}, it is evident that the center shifts from the origin to $3$, and the radius expands to approximately twice its previous size for both $d = 100$ and $d=1000$. Since~\Cref{fig: nilpotent-jordan} and~\Cref{fig: linear-symbol} are plotted on the same scale, a closer inspection reveals that the concentration radius slightly contracts when transitioning from the nilpotent Jordan block $J$ to the $\mathcal{R}$-Toeplitz matrix $T= T(3,2)$.    

To understand the observed phenomena, we analyze the characteristic matrix of $\mathcal{R}$-Toeplitz matrices under rank-1 sampling perturbations. For the $\mathcal{R}$-Toeplitz matrix defined in~\eqref{eqn: linear-symbol-teoplitz}, the associated characteristic matrix takes the form 
\[
\lambda I - T - \frac{\epsilon \pmb{u} \pmb{v}^{\dagger}}{\|\pmb{u}\|\|\pmb{v}\|} = a_1 \left[ \left( \frac{\lambda - a_0}{a_1} \right)I - J - \frac{\epsilon}{a_1} \cdot  \frac{ \pmb{u} \pmb{v}^{\dagger}}{\|\pmb{u}\|\|\pmb{v}\|}  \right].      
\]
Since both $\pmb{u}, \pmb{v} \sim \mathcal{CN}(\pmb{0},I)$ are rotationally invariant, this expression can be further simplified as
\begin{equation}
\lambda I - T - \frac{\epsilon \pmb{u} \pmb{v}^{\dagger}}{\|\pmb{u}\|\|\pmb{v}\|}  = a_1 \left[ \left( \frac{\lambda - a_0}{a_1} \right)I - J - \frac{\epsilon}{|a_1|} \cdot  \frac{ \pmb{u} \pmb{v}^{\dagger}}{\|\pmb{u}\|\|\pmb{v}\|}  \right].       \label{eqn: linear-symbol-transform}
\end{equation}
By substituting $\lambda \rightarrow \frac{\lambda - a_0}{a_1}$ and $\epsilon \rightarrow \frac{\epsilon}{|a_1|}$ into~\eqref{eqn: linear-symbol-transform}, and applying~\Cref{thm: nilpotent-jordan}, we obtain a complete mathematical characterization of the singular pseudospectral concentration for $\mathcal{R}$-Toeplitz matrices with linear symbols of the form~\eqref{eqn: linear-symbol-teoplitz}. Moreover, by performing the change of variables $\lambda \rightarrow \frac{\lambda - a_0}{a_1}$, we find that the falling set associated with the linear symbol corresponds to the image of the annuls $\left\{ \lambda  \in \mathbb{C} : \left | \left| z \right| - 1 \right | < \delta \right\} $ under the map $\lambda = p(z)$, that is,
\begin{equation}
\label{eqn: linear-symbol-tran}
\big\{ \lambda = p(z) \in \mathbb{C} : \big |  | z |  - 1\big | < \delta  \big\}  = \left\{ \lambda  \in \mathbb{C} : \left | \left|  \frac{\lambda - a_0}{a_1} \right| - 1 \right | < \delta \right\}.  
\end{equation}
For a general binomial symbol of the form~\eqref{eqn: binomial} with $m\neq 1$, the proof techniques developed in~\Cref{subsec: singular-concentration} can be extended by replacing $J$ with $J^m$ and considering the full dimension to be $md$ instead of $d$. As a result, it is straightforward to establish the corresponding result for $\mathcal{R}$-Teoplitz matrices with general binomial symbols. The result is summarized in the following theorem.

\begin{theorem}[Singular Spectral Concentration, Binomial Symbol]
\label{thm: binomial-symbol}
Let $T \in \mathscr{T}(R, d)$ be an $\mathcal{R}$-Teoplitz matrix generated by the binomial symbol~\eqref{eqn: binomial}, and let $\pmb{u}, \pmb{v} \sim \mathcal{CN}(0,1)$ two independent standard complex normal vectors. For any $\epsilon > 0$, there exists universal constants $d_0 > 0$ and $C>0$ such that, for any dimension $d \geq d_0$, the following concentration inequality holds:
\begin{equation}
\label{eqn: polynomial-symbol}
\mathbf{Pr}\left( \sigma\left( T +  \frac{  \epsilon  \pmb{u}  \pmb{v}^{\dagger}}{\|\pmb{u}\| \|\pmb{v}\|} \right)   \subseteq \big\{ \lambda = p(z) \in \mathbb{C} : \big |  |z|  - 1 \big | < \delta_d  \big\} \right)  \geq 1 - \frac{\epsilon}{|a_m|},
\end{equation}
where the concentration radius satisfies 
\[
\delta_d = 2\sqrt{\frac{C \epsilon }{  |a_m| d}}.
\]
\end{theorem}
%
%

%

\subsection{General $\mathcal{R}$-Toeplitz matrices}
\label{subsec: general}

\begin{figure}[htb!]
\centering
\begin{subfigure}[t]{0.45\linewidth}
\centering
\includegraphics[scale=0.48]{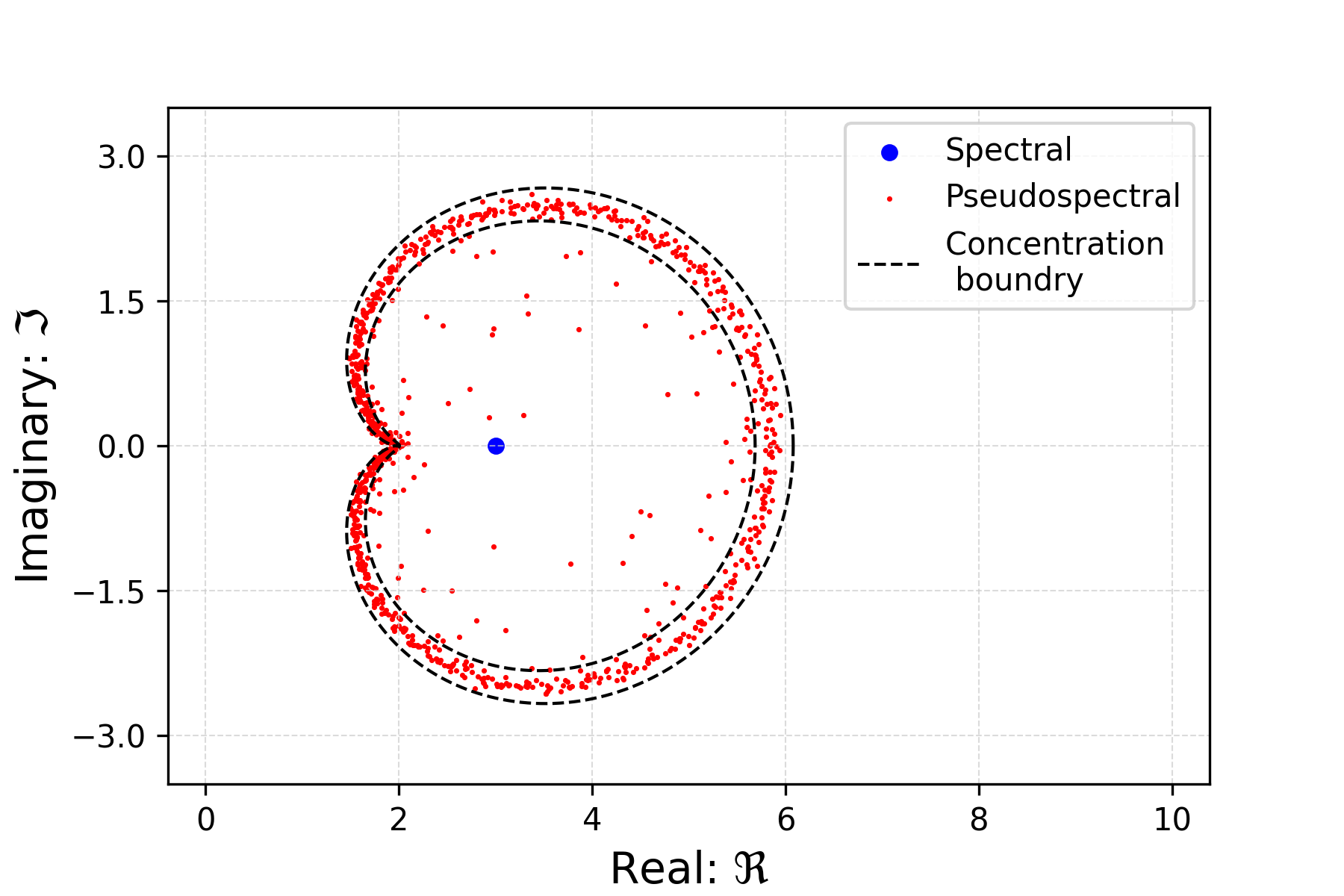}
\caption{$d=100$}
\end{subfigure}
\begin{subfigure}[t]{0.45\linewidth}
\centering
\includegraphics[scale=0.48]{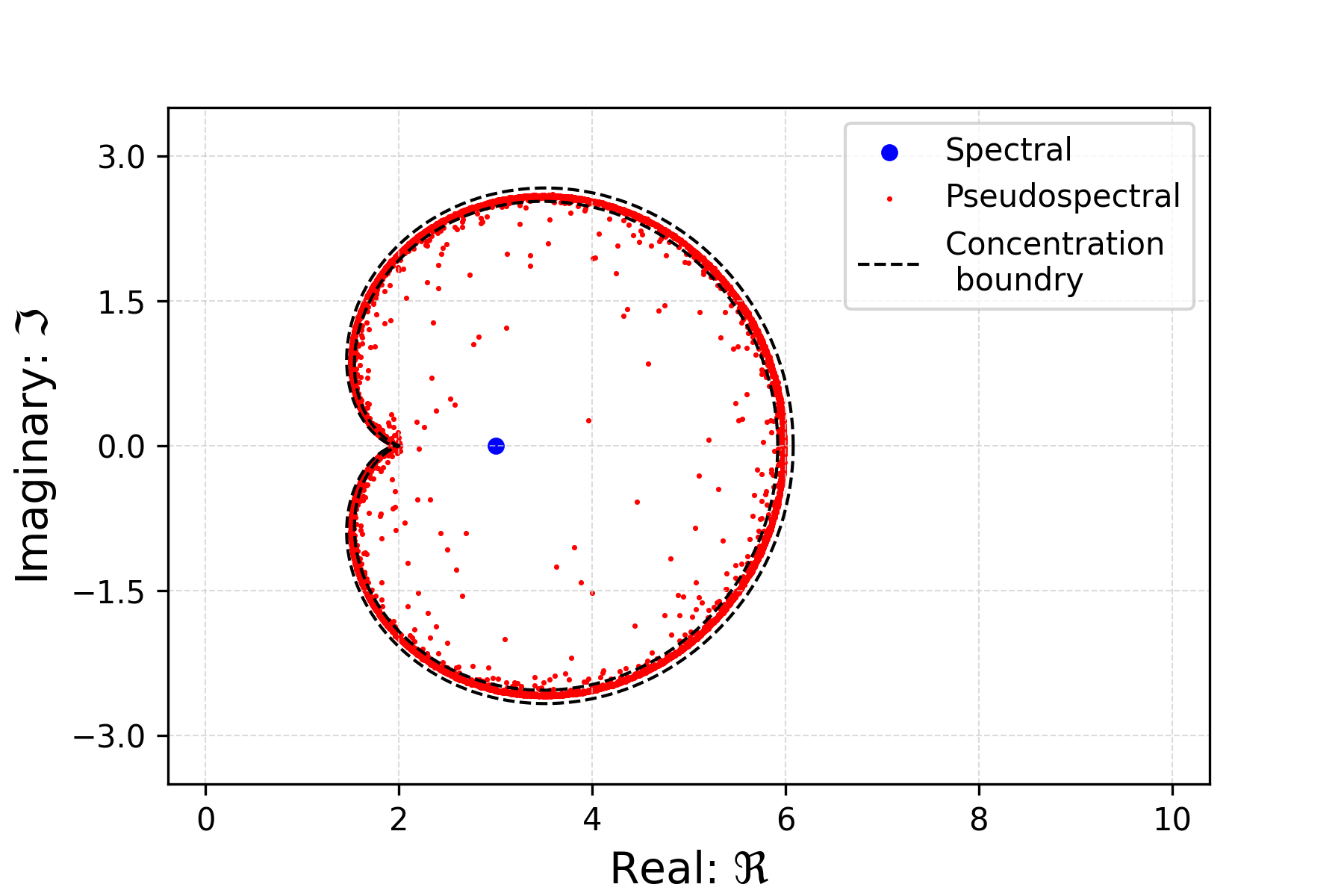}
\caption{$d=1000$}
\end{subfigure}
\caption{Scatter plots illustrating the spectral behavior of the $\mathcal{R}$-Toeplitz matrix defined by a quadratic symbol $p(z)=3+2z+z^2$ under rank-1 sampling perturbations. Parameter: $\epsilon = 2$. Number of samples: $N=1000$.} 
\label{fig: pseudospectra-polynomial}
\end{figure}

In this section, we consider $\mathcal{R}$-Toeplitz matrices generated by the general symbol~\eqref{eqn: symbol}. Unlike the binomial case, the matrix $T$ is no longer a monomial of the nilpotent Jordan block, but rather a linear combination of such monomials. As a result the resolvent $(\lambda I - T)^{-1}$ contains polynomial entries, making the analysis significantly more complex. To provide insight into the spectral behavior,  we first consider an illustrative example involving an $\mathcal{R}$-Toeplitz matrix generated by the polynomial symbol $p(z)=3+2z+z^2$ under rank-1 sampling perturbations. As depicted in~\Cref{fig: pseudospectra-polynomial}, the concentration center deviates from the circular pattern observed in the binomial case (\Cref{subsec: binomial}), instead forming a closed curve.  
This curve corresponds to the image of the unit circle  $|z|=1$ under the mapping defined by the polynomial symbol $p(z) = 3+2z+z^2$.  The shift from a simple circular distribution to a more intricate curve highlights the increased analytical challenges posed by general polynomial symbols.

To compute the resolvent, we leverage the intrinsic connection between an $\mathcal{R}$-Toeplitz matrix and its generating symbol. In this context, we introduce the characteristic equation defined by
\begin{equation}
\label{eqn: characteristic-curve}
Q(z;\lambda) = \lambda - p(z) =0,
\end{equation} 
where  $\lambda \in \mathbb{C}$ is treated as a continuously varying parameter. For certain $\lambda_0 \in \mathbb{C}$, the characteristic equation may admit multiple roots, which occur at points satisfying   
\begin{equation}
\label{eqn: multiple-root}
p'(z_{\ell}) =0, \quad \ell =1,\ldots,n-1.
\end{equation} 
In the special case of the binomial symbol $p(z) = a_0 + a_mz^m$ with $m\neq 1$, the origin is the $m$-fold roots of the symbol corresponding to $\lambda = a_0$. This singular value can be excluded by applying elementary column operations that preserve the rank, as shown in~\eqref{eqn: lambda-neq-a}. However, for a general symbol of the form~\eqref{eqn: symbol},  such multiple-root scenarios cannot be excluded so easily and must be explicitly addressed. For any $\epsilon>0$, we define the exclusion set 
\begin{equation}
\label{eqn: dig-set}
S_p(\tau) : =\bigcup_{\ell=1}^{n-1} O(p(z_{\ell}); \epsilon) = \bigcup_{\ell=1}^{n-1} \left\{ \lambda \in \mathbb{C}: |\lambda - p(z_{\ell})| < \epsilon \right\}, 
\end{equation}
which consists of neighborhoods around the critical values $p(z_\ell)$ corresponding to multiple roots of the characteristic equation. These regions must be excluded from the analysis to ensure that the modulus of the difference between any two distinct roots of the characteristic equation~\eqref{eqn: characteristic-curve} is bounded away from zero. This requirement is formalized in the following result.
\begin{lemma}
\label{lem: roots-seperation}
Let $\gamma_{\ell}(\lambda)$, for $\ell =1,\ldots,n$, denote the roots of the characteristic equation~\eqref{eqn: characteristic-curve}, and let $\epsilon>0$ be a fixed constant. Then, for all $\lambda \in \mathbb{C} \setminus S_p(\epsilon)$ (as defined in~\eqref{eqn: dig-set}), the following properties hold:
\begin{itemize}
\item[(1)] The derivatives at the roots are uniformly bounded away from zero:
                \begin{equation}
                \label{eqn: derivative-bound}
                |p'(\gamma_{\ell}(\lambda))| \geq b(\epsilon),  \qquad \mathrm{for}\;\ell = 1, \ldots, n,
                \end{equation}
                where the constant $b(\epsilon) > 0$ depends only on $\epsilon$ and is independent of $\lambda$. Moreover, 
                \begin{equation}
                \label{eqn: derivative-bound-infty}
                 b(\epsilon) \rightarrow \infty, \quad \mathrm{as}\;\; |\lambda| \rightarrow \infty.
                \end{equation}
\item[(2)]  The roots are well-separated:      
                 \begin{equation}
                \label{eqn: root-seperation}
               \min_{\ell_1 \neq \ell_2} | \gamma_{\ell_1}(\lambda) - \gamma_{\ell_2}(\lambda)| \geq \kappa(\epsilon),  \qquad \mathrm{for}\; \ell_1, \ell_2 = 1, \ldots, n,
                \end{equation}    
                 where the constant $\kappa(\epsilon) > 0$ depends only on $\epsilon$ and is independent of $\lambda$. Furthermore,
                 \begin{equation}
                \label{eqn: root-seperation-infty}
                 \kappa(\epsilon) \rightarrow \infty, \quad \mathrm{as}\;\; |\lambda| \rightarrow \infty.
                \end{equation}     
\end{itemize}
\end{lemma}

The complete proof is provided in~\Cref{sec: app-3}. Moreover, we can express the characteristic equation~\eqref{eqn: characteristic-curve} in terms of its roots as
\begin{equation}
\label{eqn: characteristic-curve-root}
Q(z; \lambda) = \lambda - p(z) = - a_n \prod_{\ell =1}^{n} \left(z - \gamma_{\ell}(\lambda) \right). 
\end{equation}
As a rational function in the form of~\eqref{eqn: characteristic-curve-root}, we can derive its partial fraction decomposition:  
\begin{equation}
\label{eqn: partial-fraction}
\frac{1}{\lambda - p(z)} = - \frac{1}{a_n} \prod_{\ell =1}^{n} \frac{1}{z - \gamma_{\ell}(\lambda) } = \sum_{\ell=1}^{n} \frac{1}{p'(\gamma_{\ell}(\lambda))(\gamma_{\ell}(\lambda) - z)}. 
\end{equation}
Since the $\mathcal{R}$-Toeplitz matrix $T$ can be expressed as the symbol polynomial evaluated at the nilpotent Jordan block $J$, we utilize~\eqref{eqn: characteristic-curve-root} to write its characteristic matrix as
\begin{equation}
\label{eqn: chracteristic-matrix-root}
\lambda I - T = \lambda I - p(J) =  - a_n \prod_{\ell=1}^{n} \left(J - \gamma_{\ell}(\lambda)I \right).
\end{equation}
Using~\Cref{lem: inverse-jordan} and the partial fraction decomposition~\eqref{eqn: chracteristic-matrix-root}, the resolvent admit the following expression: 
\begin{equation}
\label{eqn: resolvent-toeplitz-chracteristic-matrix-root}
(\lambda I - T)^{-1} =    \sum_{\ell=1}^{n} \frac{(\gamma_{\ell}(\lambda)I - J)^{-1}}{p'(\gamma_{\ell}(\lambda))}  = \sum_{k=0}^{d-1} \left( \sum_{\ell=1}^{n} \frac{1}{p'(\gamma_{\ell}(\lambda)) \gamma_{\ell}^{k+1}(\lambda) } \right) J^{k}. 
\end{equation}
Substituting the resolvent formula~\eqref{eqn: resolvent-toeplitz-chracteristic-matrix-root} into the rank-one perturbation formula, we obtain a closed-form expression for the characteristic polynomial of the perturbed matrix:
\begin{align}
 \det\left( \lambda I - T - \frac{\epsilon \pmb{u} \pmb{v}^{\dagger}}{\| \pmb{u} \| \| \pmb{v} \|} \right)    & = \det \left( \lambda I - T \right) \cdot  \left(1  -  \frac{\epsilon \pmb{v}^{\dagger} (\lambda I - T)^{-1} \pmb{u}  }{\| \pmb{v} \| \| \pmb{u} \| }  \right) \nonumber \\
                                                                                                                                                             & = (\lambda - a_0)^{d} \left( 1 - \epsilon\sum_{k=0}^{d-1} \left( \sum_{\ell=1}^{n} \frac{1}{p'(\gamma_{\ell}(\lambda)) \gamma_{\ell}^{k+1}(\lambda) } \right)  \frac{ \pmb{v}^{\dagger} J^{\ell} \pmb{u}  }{\| \pmb{v} \| \| \pmb{u} \| }  \right).      \label{eqn: close-characeteristic-toeplitz-r}
\end{align}
Similar to the analysis in~\eqref{eqn: lambda-neq-a}, $\lambda = a_0$ is also almost surely not an eigenvalue of the perturbed matrix. With the above derivations, we are now in a position to rigorously establish the singular spectral concentration behavior for the $\mathcal{R}$-Toeplitz matrix generated by the polynomial symbol~\eqref{eqn: symbol}, as stated below. 

\begin{theorem}[Singular Spectral Concentration, General Polynomial Symbol]
\label{thm: general-toeplitz}
Let $T\in \mathscr{T}(R, d)$ be an $\mathcal{R}$-Toelitz matrix generated by the polynomial symbol~\eqref{eqn: symbol}, and let $\pmb{u}, \pmb{v} \sim \mathcal{CN}(0,1)$ two independent standard complex normal vectors. For any $\epsilon > 0$, there exists universal constants $d_0 > 0$ and $C>0$ such that, for any dimension $d \geq d_0$, the following concentration inequality holds:
\begin{equation}
\label{eqn: polynomial-symbol}
\mathbf{Pr}\left( \sigma\left( T +  \frac{  \epsilon  \pmb{u}  \pmb{v}^{\dagger}}{\|\pmb{u}\| \|\pmb{v}\|} \right)   \subseteq \big\{ \lambda = p(z) \in \mathbb{C} : \big |  |z|  - 1\big | < \delta_d \big\} \cup S_p(\epsilon) \right)  \geq 1 - \frac{\epsilon}{|a_n|},
\end{equation}
where the concentration radius satisfies 
\[
\delta_d = 2\sqrt{\frac{ C \epsilon }{ nb(\tau_0)|a_n|d}}.
\] 
\end{theorem}

\begin{proof}[Proof of~\Cref{thm: general-toeplitz}]
We establish the probability that all eigenvalues lie within the annular region $\big\{ \lambda = p(z) \in \mathbb{C} : \big |  | z |  - 1\big | < \delta  \big\}$ through the decomposition of the characteristic roots. The probability can be written as 
\[
\mathbf{Pr} \bigg( \big\{ \lambda = p(z) \in \mathbb{C} : \big |  | z |  - 1\big | < \delta  \big\}  \bigg)     = \mathbf{Pr} \left( \bigcup_{\ell =1}^{n} \left\{ \gamma_k(\lambda) \in \mathbb{C}:  \left| |\gamma_{\ell}(\lambda)| - 1\right| < \delta\right\}   \right). 
\]
Furthermore, this expression admits a refined decomposition into two dominant components: 
\begin{multline}
 \mathbf{Pr} \bigg( \big\{ \lambda = p(z) \in \mathbb{C} : \big |  | z |  - 1\big | < \delta  \big\}  \bigg)                                                                      \\
 =  \underbrace{\mathbf{Pr} \left( \bigcup_{\ell=1}^{n} \left\{ \gamma_k(\lambda) \in \mathbb{C}:  |\gamma_{\ell}(\lambda)| < 1+ \delta\right\} \right)}_{:=\mathbf{I}}  - \underbrace{\mathbf{Pr} \left( \bigcap_{\ell=1}^{n} \left\{ \gamma_{\ell}(\lambda) \in \mathbb{C}:  |\gamma_{\ell}(\lambda)| \leq 1- \delta\right\}  \right)}_{:=\mathbf{II}}.      \label{eqn: r-toeplitz-separation}
\end{multline}
where the first term $\mathbf{I}$  captures configurations where at least one characteristic root penetrates the outer boundary $1+\delta$, and the second term $\mathbf{II}$ involves the case where all roots remain confined within the inner core $1-\delta$. 


\begin{itemize}
\item[(1)] \textbf{At least one root penetrates the outer boundary.} According to the characteristic polynomial in~\eqref{eqn: close-characeteristic-toeplitz-r}, if $\lambda$ is a root, then it must satisfy the following condition:
\begin{equation}
\label{eqn: general-polynomial-outer-1}
1 - \epsilon\sum_{k=0}^{d-1} \left( \sum_{\ell=1}^{n} \frac{1}{p'(\gamma_{\ell}) \gamma_{\ell}^{k+1} } \right)  \frac{ \pmb{v}^{\dagger} J^{\ell} \pmb{u}  }{\| \pmb{v} \| \| \pmb{u} \| }  = 0,
\end{equation}
where, for convenience, we suppress the dependence on $\lambda$ and write $\gamma_{\ell}$ instead of $\gamma_{\ell}(\lambda)$.  According to~\Cref{lem: maximum-modulus-principle} (Maximum Modulus Principle), in order for the event included in $\mathbf{I}$ to be active, at least one root must lie on the outer boundary, i.e., $|\gamma_{\ell}| = 1 + \delta$. Under this condition, the left-hand side of~\eqref{eqn: general-polynomial-outer-1} must be strictly less than one:
\[
 \epsilon \left| \sum_{k=0}^{d-1} \left( \sum_{\ell=1}^{n} \frac{1}{p'(\gamma_{\ell}) \gamma_{\ell}^{k+1} } \right)  \frac{ \pmb{v}^{\dagger} J^{\ell} \pmb{u}  }{\| \pmb{v} \| \| \pmb{u} \| }  \right| < 1. 
 \]
Conversely, if no such root exists, then the complement case must satisfy
\[
 \epsilon \left| \sum_{k=0}^{d-1} \left( \sum_{\ell=1}^{n} \frac{1}{p'(\gamma_{\ell}) \gamma_{\ell}^{k+1} } \right)  \frac{ \pmb{v}^{\dagger} J^{\ell} \pmb{u}  }{\| \pmb{v} \| \| \pmb{u} \| }  \right| \geq 1.
 \]
Applying the triangle inequality, we can bound the probability $\mathbf{I}$ from below as
\begin{equation}
\label{eqn: general-polynomial-outer-2}
\mathbf{I} \geq 1- \mathbf{Pr} \left( \epsilon  \sum_{k=0}^{d-1} \left| \mathbf{v}^{\dagger} \Omega_k \pmb{u}\right| \geq  \| \pmb{v} \| \| \pmb{u} \|  \right).
 \end{equation}
where the perturbation terms are defined as
\[
\Omega_k := \frac{ J^k  }{(1+\delta)^{k+1}} \left( \sum_{\ell=1}^{n} \frac{1}{p'(\gamma_\ell) e^{i(k+1)\theta_\ell}} \right), \quad k = 0,1,\ldots, d-1.
\]
Let $\tau_0>0$ be a parameter to be determined later.  For $\lambda \in \mathbb{C} \setminus  S_p(\tau_0)$, we further relax the lower bound
\begin{equation}
\label{eqn: general-polynomial-outer-3}
\mathbf{I} \geq 1- \mathbf{Pr} \left( \frac{\epsilon}{nb(\tau_0)}   \sum_{k=0}^{d-1} \left| \frac{\mathbf{v}^{\dagger}J^k \pmb{u}}{(1+\delta)^{k+1}} \right| \geq  \| \pmb{v} \| \| \pmb{u} \|  \right).
\end{equation}
The remaining analysis follows identically to the derivation from equations~\eqref{eqn: linear-symbol-outer-2} to~\eqref{eqn: linear-symbol-outer-8}, with the only modification being the use of $\frac{\epsilon}{nb(\tau_0)} $ in place of $\epsilon$. In particular, the separation radius is given by:
\[
\delta_d = 2\sqrt{\frac{ C \epsilon }{ nb(\tau_0)|a_n|d}}.
\]

\item[(2)] \textbf{All roots remain confined within the inner core.}  We now aim to demonstrate that the event where all roots remain confined within the inner core occurs with small probability. To proceed, we invoke~\Cref{lem: roots-seperation}, which asserts that at most one of the roots $\gamma_k$ can be zero. This guarantees the existence of a strictly positive minimal modulus among the roots, which is defined as
\begin{equation}
\label{eqn: root-min-teoplitz}
\mu := \min_{1 \leq \ell \leq n}|\gamma_{\ell}|.
\end{equation}
According to~\Cref{lem: maximum-modulus-principle} (Maximum Modulus Principle), we estimate the probability of this inward deviation by establishing a bound on the boundary. Specifically, we have 
\begin{equation}
\label{eqn: general-polynomial-inner-1}
\mathbf{II} \leq \mathbf{Pr} \left( \epsilon  \left| \sum_{k=0}^{d-1}  \pmb{v}^{\dagger} \omega_k \pmb{u} \right| \leq  \| \pmb{v} \| \| \pmb{u} \|  \right)
 \end{equation}
where the perturbation terms are given by
\[
\omega_k :=   \left( \sum_{\ell=1}^{n} \frac{1}{p'(\gamma_\ell) \gamma_{\ell}^{k+1} } \right)J^k, \quad k = 0,1,\ldots, d-1.
\] 

The subsequent analysis proceeds in the same manner as in~\eqref{eqn: linear-symbol-inner-1} ---~\eqref{eqn: linear-symbol-inner-5}. The only additional challenge lies in estimating a lower bound on the Frobenius norm, which is required to apply~\Cref{lem: carbery-wright-complex}. For convenience, let us define $k_1 := \left\lfloor  \frac{d}{2} \right\rfloor $, and proceed to estimate the Frobenius norm as follows:
\begin{align}
\left\| \sum_{k=0}^{d-1} \omega_k \right\|_{F}^{2} & = \sum_{k=0}^{d-1} (d-k) \left| \sum_{\ell=1}^{n} \frac{1}{p'(\gamma_\ell) \gamma_{\ell}^{k+1} } \right|^2 \nonumber \\
                                                                              & \geq  \sum_{k=0}^{d-1} \frac{1}{\mu^{2(k+1)}} \left| \sum_{\ell=1}^{n} \frac{\mu^{k+1}}{p'(\gamma_\ell)  \gamma_{\ell}^{k+1} } \right|^2 \geq  \frac{1}{\mu^{d}} \sum_{k=k_1}^{d-1} \left| \sum_{\ell=1}^{n} \frac{\mu^{k+1}}{p'(\gamma_\ell)  \gamma_{\ell}^{k+1} } \right|^2.  \label{eqn: general-polynomial-inner-2}
\end{align}
Furthermore, by expanding the squared modulus more carefully, we obtain:
\begin{align}
\left\| \sum_{k=0}^{d-1} \omega_k \right\|_{F}^{2}  \geq & \underbrace{\frac{1}{\mu^{d}} \sum_{k=k_1}^{d-1} \sum_{\ell=1}^{n} \frac{\mu^{2(k+1)}}{|p'(\gamma_{\ell})|^2 |\gamma_{\ell}|^{2(k+1)}}}_{:=\mathbf{K}_1} \nonumber \\ 
&- \underbrace{\frac{1}{\mu^d} \left|  \sum_{k=k_1}^{d-1} \sum_{\substack{\ell_1 \neq \ell_2 \\ 1 \leq \ell_1, \ell_2 \leq n}} \frac{\mu^{2(k+1)}}{p'(\gamma_{\ell_1})\overline{p'(\gamma_{\ell_2})} \left( \gamma_{\ell_1}  \overline{\gamma}_{\ell_2} \right)^{k+1}}  \right|}_{:=\mathbf{K}_2}.  \label{eqn: general-polynomial-inner-3}
\end{align}
Since all roots $\gamma_k$ lie within the unit disk $\mathbb{D} = \{ z \in \mathbb{C}: |z|\leq1\}$, the derivatives $p'(\gamma_k)$ are uniformly bounded on $\mathbb{D}$. This bound is given by the maximum modulus of $p'$ on the disk:
\begin{equation}
\label{eqn: pderi-max-teoplitz}
L = \max_{|z| \leq 1} |p'(z)|. 
\end{equation}
For the first term $\mathbf{K}_1$ in~\eqref{eqn: general-polynomial-inner-3}, using the definition of the minimum modulus $\mu$ from~\eqref{eqn: root-min-teoplitz} and the maximum derivative bound $L$ from~\eqref{eqn: pderi-max-teoplitz}, we can derive the following lower bound: 
\begin{equation}
\label{eqn: k1-lower-toeplitz}
K_1 \geq \frac{1}{\mu^d} \sum_{k=k_1}^{d-1} \frac{1}{L} \geq  \frac{nd}{2\mu^dL}. 
\end{equation}
Consider the second term $\mathbf{K}_1$ in~\eqref{eqn: general-polynomial-inner-3}, we  swap the order of summation to further estimate the upper bound as: 
\begin{align}
K_2 & = \frac{1}{\mu^d} \left| \sum_{\substack{\ell_1 \neq \ell_2 \\ 1 \leq \ell_1, \ell_2 \leq n}} \frac{\mu^{2(k_1+1)}}{p'(\gamma_{\ell_1})\overline{p'(\gamma_{\ell_2})} ( \gamma_{\ell_1} \overline{\gamma}_{\ell_2} )^{2(k_1+1)} } \left[ \sum_{k=0}^{d-k_1-1}\left( \frac{\mu^2}{\gamma_{\ell_1} \overline{\gamma}_{\ell_2}} \right)^k \right] \right|  \nonumber \\
       & \leq \frac{2}{\mu^d} \sum_{\substack{\ell_1 \neq \ell_2 \\ 1 \leq \ell_1, \ell_2 \leq n}} \frac{1}{|p'(\gamma_1)||p'(\gamma_2)|} \cdot \frac{1}{\left|1 - \frac{\mu^2}{\gamma_{\ell_1} \overline{\gamma}_{\ell_2}} \right|},        \label{eqn: k2-upper-toeplitz} 
\end{align}
where the inequality also follows the definition of the minimum modulus $\mu$ from~\eqref{eqn: root-min-teoplitz}. If the final term in~\eqref{eqn: k2-upper-toeplitz} is uniformly bounded for all $\ell_1 \neq \ell_2$, specifically, 
\begin{equation}
\label{eqn: final-estimate-Toeplitz}
\left|1 - \frac{\mu^2}{\gamma_{\ell_1} \overline{\gamma}_{\ell_2}} \right| \geq \eta(\epsilon) \qquad \mathrm{for}\;\; \ell_1 \neq \ell_2, \;\; \ell_1,\ell_2 = 1, \ldots,n;
\end{equation}
then, by combining the two bounds~\eqref{eqn: k1-lower-toeplitz} and~\eqref{eqn: k2-upper-toeplitz}, we obtain the following lower bound on the Frobenius norm:
\begin{equation}
\label{eqn: general-polynomial-inner-4}
\left\| \sum_{k=0}^{d-1} \omega_k \right\|_{F}^{2}  \geq \frac{n}{\mu^d} \left(\frac{d}{2L} -  \frac{2n}{\eta(\epsilon) b^{2}(\epsilon)}\right).
\end{equation}
This implies that the Frobenius norm grows asymptotically at a rate no less than $\frac{d}{\mu^d}$ with $\mu \leq 1 - \delta$, thereby ensuring that the inequality~\eqref{eqn: linear-symbol-inner-5} always holds. Finally, we turn to estimating the bound in~\eqref{eqn: final-estimate-Toeplitz}. Using the identity for the square modulus, we have the expansion:
\[
2\Re\left( \frac{\mu^2}{\gamma_{\ell_1} \overline{\gamma}_{\ell_2}} \right)= \left| \frac{\mu}{\gamma_{\ell_1}} \right|^2 + \left| \frac{\mu}{\gamma_{\ell_2}} \right|^2  - \left| \frac{\mu}{\gamma_{\ell_1}} - \frac{\mu}{\gamma_{\ell_2}} \right|^2. 
\]
This allows us to estimate the final term in~\eqref{eqn: final-estimate-Toeplitz} as:
\[
\left| 1 - \frac{\mu^2}{\gamma_{\ell_1} \overline{\gamma}_{\ell_2}} \right| \geq 1 - \frac12 \left(\left| \frac{\mu}{\gamma_{\ell_1}} \right|^2 + \left| \frac{\mu}{\gamma_{\ell_2}} \right|^2\right)  + \frac{\left| \gamma_{\ell_1} -\gamma_{\ell_2} \right|^2}{2|\gamma_{\ell_1} | |\gamma_{\ell_2} |}. 
\]
We now divide the analysis into two distinct cases: first, when both $|\gamma_{\ell_1}|$ and $|\gamma_{\ell_2}|$ are less than $2\mu$, and second at least one of them is no less than $2\mu$. Across both scenarios, we establish a universal lower bound: 
\[
\left| 1 - \frac{\mu^2}{\gamma_{\ell_1} \overline{\gamma}_{\ell_2}} \right| \geq \min\left\{ \frac{3}{8}, \kappa(\epsilon)\right\}.
\]
Thus, by setting $\eta(\epsilon): =  \min\left\{ 3/8, \kappa(\epsilon)\right\}$,  we conclude  the derivation of the lower bound in~\eqref{eqn: final-estimate-Toeplitz}. Hence, the proof is complete. 
\end{itemize}
\end{proof}

\section{Conclusion and future work}
\label{sec: conclusion}

In this study, we  have established a theoretical framework to quantify the pseudospectral concentration behavior of complex matrices under rank-$1$ sampling perturbations. In particular, for singular concentration phenomena, we derived a sub-exponential tail bound for nilpotent Jordan blocks by combining the Hanson–Wright concentration inequality and the Carbery–Wright anti-concentration inequality. By analyzing root-separation properties of polynomials in the complex plane, we identified and excluded regions where arguments of distinct roots exhibit insufficient separation, thereby extending the tail bound to $\mathcal{R}$-Toeplitz matrices. Numerical experiments indicate that degenerate cases involving multiple roots rarely influence the overall singular concentration behavior. Nevertheless, the root-separation region near such degeneracies presents mathematically challenging. Within these regions, partial fraction decomposition remains applicable but requires careful handling of higher-order derivatives. Specifically, for a polynomial $\lambda = p(z)$ with roots $\{\gamma_{\ell}\}_{\ell=1}^{r}$ of multiplicities $s_{\ell}$, the decomposition takes the form:
\[
\frac{1}{\lambda - p(z)} = \sum_{\ell=1}^{r} \sum_{k=1}^{s_{\ell}} \frac{A_{\ell,k}}{(\gamma_{\ell} - z)^k}, \quad \mathrm{where}\;\;A_{\ell, k} = \frac{1}{(s_{\ell} -k)} \lim_{z \rightarrow \gamma_{\ell}} \frac{d^{s_{\ell} -k}}{dz^{s_{\ell} - k}} \left[ \frac{(\gamma_{\ell} - z)^k}{\lambda - p(z)}  \right].
\]
This formulation potentially refines the exclusion of poorly separated root regions by leveraging relationships between roots and coefficients of $p(z)$.

Numerical evidence suggests that the behavior of low-rank and even full-rank sampling perturbations is nearly consistent with the rank-1 perturbation case.  A promising direction for future work is to rigorously analyze such perturbations. Building upon tools from random matrix theory (e.g.,~\citep{tao2012topics}), recent studies such as~\citet{guionnet2014convergence} and~\citet{basak2020spectrum} have explored weak convergence in the sense of distribution. Extending this to a quantitative theory that characterizes the statistical properties of general low-rank (or full-rank) sampling perturbations would be highly valuable. Such a theory could exploit the fundamental decomposition of perturbations into normalized rank-$1$ components:
\[
T +  \epsilon  \left( \frac{  \pmb{u}_1  \pmb{v}_1^{\dagger}}{\|\pmb{u}_1 \| \|\pmb{v}_1\|} + \cdots +   \frac{\pmb{u}_m  \pmb{v}_m^{\dagger}}{\|\pmb{u}_m \| \|\pmb{v}_m\|} \right). 
\]
Complementing this theoretical direction, random sampling techniques offer a practical and efficient approach to numerically investigate pseudospectral behavior in nonnormal operators. The seminal work of~\citep{reddy1993pseudospectra} demonstrated the effectiveness of such methodsin revealing nonnormal effects in hydrodynamic stability problems. These computational insights have since inspired  rigorous analytical developments~\citep{ibrahim2019pseudospectral, li2020oseen, li2020pseudospectral}, particularly in fluid dynamical systems where nonnormality plays a fundamental role. Further exploration of this interplay between numerical experimentation and theoretical analysis promises to yield valuable insights for both mathematical theory and physical applications. 



%
%

%
%
%

\section*{Acknowledgments}
K. Gai sincerely expresses gratitude to Shihua Zhang for his support, which facilitated Gai's visit to SIMIS and Fudan University. This work was partially supported by startup fund from SIMIS and Grant No.12241105 from NSFC.

\bibliographystyle{abbrvnat}
\bibliography{ref}

\appendix
\section{Proof of~\Cref{lem: norm-concentration}}
\label{sec: app-1}
In this section, we establish a concentration inequality for the squared norm of a standard complex normal vector $\pmb{\xi} \sim \mathcal{CN}(\pmb{0}, I)$ using classical moment-generating function techniques. For any $\lambda \in (0, 1)$, applying Markov's inequality to the exponential moment of the squared norm, we  obtains:
\begin{equation}
\label{eqn: norm-concen-1}
\mathbf{Pr}\left( \frac{1}{d} \| \pmb{\xi}\|^2 - 1 \leq - t \right)  = \mathbf{Pr}\left(  \prod_{k=1}^{d}  e^{\lambda(1 - |\xi_k|^2 ) } \geq e^{\lambda dt} \right) \leq e^{-\lambda d (t-1)}\mathbb{E}\left[ \prod_{k=1}^{d} e^{ -\lambda |\xi_k|^2 } \right].
\end{equation}
Since the components $\xi_k$, $k=1,\ldots, d$, are independent and identically distributed, the expectation factorizes: 
\[
\mathbb{E}\left[ \prod_{k=1}^{d} e^{ - \lambda |\xi_k|^2 } \right] =   \prod_{i=1}^{d}  \mathbb{E}\left[e^{ - \lambda |\xi_k|^2} \right].  
\]
Now observe that for a standard complex normal variable $\xi_k \sim \mathcal{CN}(0,1)$, the squared modulus $|\xi_k|^2$ follows an exponential distribution with mean $1$. Thus, for any $\lambda \in (0,1)$, the moment-generating function is given by:
\[
 \mathbb{E}\left[e^{- \lambda |\xi_k|^2} \right] =  \int_{0}^{\infty} e^{- (1 + \lambda) r} rdr = \frac{1}{1+\lambda}. 
\]
Therefore, the full expectation evaluates to:
\begin{equation}
\label{eqn: norm-concen-2}
\mathbb{E}\left[ \prod_{i=1}^{d} e^{ - \lambda |\xi_k|^2 } \right] =  \frac{1}{(1 + \lambda)^d}. 
\end{equation}
Substituting~\eqref{eqn: norm-concen-2} into the inequality~\eqref{eqn: norm-concen-1}, we obtain the lower-tail bound
\[
\mathbf{Pr}\left( \frac{1}{d} \| \pmb{\xi}\|^2 - 1 \leq -t \right)  \leq  \exp\left(-\lambda d (t-1) - d \log (1 + \lambda )\right).
\]
To simplify this expression, we invoke the inequality $\log(1 + \lambda) \geq \lambda - \lambda^2/2 $ for any $\lambda \in (0,1)$, which yields:
\begin{equation}
\label{eqn: norm-concen-4}
\mathbf{Pr}\left( \frac{1}{d} \| \pmb{\xi}\|^2 - 1\leq -t \right)  \leq \exp\left(-\lambda d t + \frac{\lambda^2 d}{2} \right).
\end{equation}
By symmetry, a similar argument gives the upper-tail bound: 
\begin{equation}
\label{eqn: norm-concen-5}
\mathbf{Pr}\left( \frac{1}{d} \| \pmb{\xi}\|^2 - 1 \geq t \right)  \leq \exp\left(-\lambda d t  \right),
\end{equation}
which also holds for any $\lambda \in (0, 1)$.  Finally, setting $\lambda = \sqrt{1/d}$ in both bounds~\eqref{eqn: norm-concen-4} and~\eqref{eqn: norm-concen-5}, we complete the proof with a universal constant $C =1+\sqrt{e}$.

\section{Technical details in~\Cref{subsec: complex-probabilstic-ineq}}
\label{sec: technical-detail-complex-ineq}

In this section, we provide detailed proofs of the two key probabilistic inequalities in the complex setting, as introduced in~\Cref{subsec: complex-probabilstic-ineq}:  the Hanson-Wright concentration inequality~(\Cref{lem: hanson-wright-complex}) and the Carbery–Wright anti-concentration inequality~(\Cref{lem: carbery-wright-complex}).



\subsection{Proof of~\Cref{lem: hanson-wright-complex}} 
\label{subsec: app-21}

This proof begins by analyzing the complex quadratic forms $\pmb{v}^{\dagger}J_0^k \pmb{u}$ for $k=0,1,\ldots,d-1$. By decomposing each variable into its  real and imaginary parts, we write $v_{\ell} = v_{\ell, r} + i v_{\ell, i}$ and $u_{k+\ell} = u_{k+\ell, r} + i u_{k+\ell,i}$, which leads to the following expressions for the real and imaginary parts of the quadratic forms: 
\begin{subequations}
\begin{align}
&  \Re \left( \pmb{v}^{\dagger}J_0^k \pmb{u} \right) =  \sum_{\ell =1}^{d-k}    \left( v_{\ell, r}u_{k+\ell, r} + v_{\ell,i}  u_{k+\ell,i} \right),             \label{eqn: hw-real}\\
&  \Im \left( \pmb{v}^{\dagger}J_0^k \pmb{u} \right) =   \sum_{\ell =1}^{d-k}    \left( v_{\ell,r} u_{k+\ell,i}  -  v_{\ell,i}  u_{k+\ell,r} \right).            \label{eqn: hw-imag}
\end{align}
\end{subequations}
Given that both $\pmb{u}$ and $\pmb{v}$ are independent standard complex normal vectors, i.e., $\pmb{u}, \pmb{v} \sim  \mathcal{CN}(\pmb{0},I)$, it follows that their real and imaginary components $u_{k,r}, u_{k,i}, v_{k,r}, v_{k,i} \sim  \mathcal{N}(0,1/2)$ are i.i.d.~real normal random variables. This structure ensures, by linearity and symmetry, that both the real and imaginary parts of $\pmb{v}^{\dagger}J_0^k \pmb{u}$ inherit zero mean:
\[
\mathbb{E}\left[  \Re \left( \pmb{v}^{\dagger}J_0^k \pmb{u}\right) \right] = \mathbb{E}\left[  \Im \left( \pmb{v}^{\dagger}J_0^k \pmb{u} \right) \right] = 0.
\]
To establish concentration bounds, we construct a real-valued random vector as
\[
\pmb{X} = \sqrt{2} (v_{1, r}, v_{1, i}, \ldots, v_{d-k, r}, v_{d-k, i}, u_{k+1, r}, u_{k+1, i}, \ldots, u_{d, r}, u_{d, i})^{\top}  \in \mathbb{R}^{4(d-k)},
\] 
which follows the standard real normal distribution, i.e., $\pmb{X} \sim \mathcal{N}(\pmb{0}, I)$.  The analysis proceeds by defining a matrix $A \in \mathscr{C}(4(d-k))$, consisting of repeated blocks, as
\[
A = \begin{array}{c}
      \begin{pmatrix}
      \begin{array}{ccc:ccc}
         \ddots     &             &             &             &               &           \\
                        & \ddots  &              &             &              &           \\
                        &             & \ddots   &             &              &           \\     \hdashline    
                  B    &             &             & \ddots   &              &          \\
                        & \ddots   &             &             &  \ddots   &          \\
                        &              &   B        &             &              & \ddots          
       \end{array}
       \end{pmatrix} \\
       \underbrace{\hspace{2.5cm}}_{2(d-k)} \underbrace{\hspace{2.5cm}}_{2(d-k)}
       \end{array} \quad \mathrm{where}\;\; B = \frac{1}{2} \begin{pmatrix}
              1 & 0               \\
              0 & 1   
\end{pmatrix}.
\]
This construction allows us to elegantly express the real part $\Re \left( \pmb{v}^{\dagger}J_0^k \pmb{u} \right)$ from~\eqref{eqn: hw-real} as a quadratic form: \begin{equation}
\label{eqn: hw-real-equal}
\Re \left( \pmb{v}^{\dagger}J_0^k \pmb{u} \right) = \pmb{X}^{\top}A \pmb{X}. 
\end{equation}
By applying~\Cref{lem: hanson-wright} (Hanson-Wright Concentration Inequality), we know that there exists a universal constant $C>0$ such that for any $t>0$, it holds for the following tail bound: 
\begin{equation}
\label{eqn: hw-real-final}
\mathbf{Pr}\left( \left| \Re \left( \pmb{v}^{\dagger}J_0^k \pmb{u}\right) \right| \geq \frac{t}{\sqrt{2}} \right) \leq 2 \exp \left( -C \min \left( \frac{t^2}{d-k}, t \right)\right).
\end{equation}
A similar argument applies to the imaginary part $\Im(\pmb{v}^{\dagger}J_0^k \pmb{u})$ from~\eqref{eqn: hw-imag}, yieldling a parallel tail bound:
\begin{equation}
\label{eqn: hw-imag-final}
\mathbf{Pr}\left( \left| \Im \left( \pmb{v}^{\dagger}J_0^k \pmb{u}\right) \right| \geq \frac{t}{\sqrt{2}} \right) \leq 2 \exp \left( -C \min \left( \frac{t^2}{d-k}, t \right)\right),
\end{equation}
Combining the two probabilistic inequalities,~\eqref{eqn: hw-real-final} and~\eqref{eqn: hw-imag-final}, we complete the proof. 

%
%

\subsection{Proof of~\Cref{lem: carbery-wright-complex}} 
\label{subsec: app-22}

As demonstrated in~\Cref{subsec: app-21}, we analyze the quadratic form $\pmb{v}^{\dagger} A \pmb{u}$ by decomposing each complex variable into its real and imaginary parts: $v_{\ell} = v_{\ell, r} + i v_{\ell, i}$ and $u_{k+\ell} = u_{k+\ell, r} + i u_{k+\ell,i}$. This decomposition allows us to express the real and imaginary parts of the quadratic form explicitly, leading to the following expansions:
\begin{subequations}
\begin{align}
&  \Re \big( \pmb{v}^{\dagger} A \pmb{u} \big) =   \sum_{k=1}^{d} \sum_{\ell=1}^{d}  \left[ a_{k\ell,r} \left( v_{k, r}u_{\ell, r}+ v_{k,i}u_{\ell,i} \right) - a_{k\ell, i}\left( v_{k,r} u_{\ell, i} -  v_{k,i}  u_{\ell,r}\right) \right],                        \label{eqn: cw-real}        \\
&  \Im \big( \pmb{v}^{\dagger} A \pmb{u} \big) =  \sum_{k=1}^{d} \sum_{\ell=1}^{d}\left[ a_{k\ell,r}\left( v_{k,r} u_{\ell,i} -  v_{k,i}  u_{\ell,r} \right) + a_{k\ell,i}\left( v_{k,r}u_{\ell,r} + v_{k,i}u_{\ell,i} \right) \right].                        \label{eqn: cw-imag}
\end{align}
\end{subequations}
Since $\pmb{u}, \pmb{v} \sim  \mathcal{CN}(\pmb{0},I)$ are independent standard complex normal vectors, their real and imaginary parts  $u_{k,r}, u_{k,i}, v_{k,r}, v_{k,i} \sim  \mathcal{N}(0,1/2)$ i.i.d. real normal random variables. Thus, the expectations relevant for computing the variance are:
\begin{subequations}
\begin{align}
&  \omega_{k\ell,r} = \mathbb{E}\left[ v_{k,r}^2 \right] \mathbb{E}\left[ u_{\ell, r}^2  \right] + \mathbb{E}\left[v_{k,i}^2 \right] \mathbb{E}\left[ u_{\ell, i}^2 \right] = \frac12,  \label{eqn: cw-omega-real}     \\
&  \omega_{k\ell,i} = \mathbb{E}\left[ v_{k,r}^2 \right] \mathbb{E}\left[  u_{\ell, i}^2  \right] + \mathbb{E}\left[v_{k,i}^2 \right] \mathbb{E}\left[ u_{\ell,r}^2 \right] = \frac12.  \label{eqn: cw-omega-imag} 
\end{align}
\end{subequations}
By symmetry and the linearity of expectation, both the real and imaginary parts of the quadratic form have zero mean: 
\begin{equation}
\label{eqn: cw-mean}
\mathbb{E}\left[  \Re \big( \pmb{v}^{\dagger} A \pmb{u} \big) \right] = \mathbb{E}\left[  \Im \big( \pmb{v}^{\dagger} A \pmb{u} \big) \right] = 0.
\end{equation}
Using the identities~\eqref{eqn: cw-omega-real},~\eqref{eqn: cw-omega-imag}, and~\eqref{eqn: cw-mean}, we compute the variance of the real part~\eqref{eqn: cw-real} as:
\begin{align}
 \mathrm{Var}\left[  \Re \big( \pmb{v}^{\dagger} A \pmb{u} \big)  \right]  & = \mathbb{E}\left[  \Re \big( \pmb{v}^{\dagger} A \pmb{u} \big)^2  \right]          \nonumber   \\
                                                                                                               &  = \sum_{k=1}^{d} \sum_{\ell=1}^{d} \bigg( a_{k\ell,r}^2 \omega_{k\ell,r} +  a_{k\ell, i}^2  \omega_{k\ell,i}  \bigg)   = \frac{1}{2} \sum_{k=1}^{d} \sum_{\ell=1}^{d} |a_{k\ell}|^2 = \frac{1}{2} \|A\|_{F}^{2} .                                                                                                           \label{eqn: cw-variance}
\end{align}
By applying~\Cref{lem: carbery-wright} (Carbery–Wright Anti-Concentration Inequality), we obtain a tail bound for the real part~\eqref{eqn: cw-real} :
\begin{equation}
\label{eqn: cw-real-final}
\mathbf{Pr}\left( \left| \Re \big( \pmb{v}^{\dagger} A \pmb{u} \big) \right| \leq \frac{t}{\sqrt{2}} \|A\|_{F}\right) \leq 2^{\frac14}C_1\sqrt{t},
\end{equation}
for some universal constant $C_1>0$. An identical argument yields the same bound for the imaginary part~\eqref{eqn: cw-imag}:
\begin{equation}
\label{eqn: cw-imag-final}
\mathbf{Pr}\left( \left| \Im \big( \pmb{v}^{\dagger} A \pmb{u} \big) \right| \leq \frac{t}{\sqrt{2}} \|A\|_{F}\right) \leq 2^{\frac14}C_1\sqrt{t}.
\end{equation}
Finally, combining the two inequalities,~\eqref{eqn: hw-real-final} and~\eqref{eqn: hw-imag-final}, and using the identity $|\pmb{v}^{\dagger} A \pmb{u}|^2 =   \Re \big( \pmb{v}^{\dagger} A \pmb{u} \big)^2 +  \Im \big( \pmb{v}^{\dagger} A \pmb{u} \big)^2$, we complete the proof by setting the universal constant $C=2^{\frac{5}{4}}C_1$.

\section{Proof of~\Cref{lem: roots-seperation}}
\label{sec: app-3}

Let $p(z)$ be the degree-$n$ polynomial defined in~\eqref{eqn: symbol}, and let $\gamma_{\ell}(\lambda)$, for $\ell =1, \ldots, n$, denote the $n$ roots of the characteristic equation $Q(z; \lambda) =0$ given in~\eqref{eqn: characteristic-curve}. We divide the proof into four steps.

\subsection*{Step 1: Analyticity and asymptotic growth of roots}
By~\Cref{lem: rouche} (Rouch\'e's Theorem), the roots $\gamma_{\ell}(\lambda)$, for $\ell = 1,\ldots,n$, depend continuously on $\lambda \in \mathbb{C}$. In fact, since $Q(z; \lambda)$ is a polynomial in $z$ with coefficients that are analytic in $\lambda$, the implicit function theorem ensures that $\gamma_{\ell}(\lambda)$ is holomorphic in $\lambda$, except possibly at branching points. These branching points can be treated systematically by considering the roots as multivalued functions on an appropriate Riemann surface. For our purposes, however, it suffices to consider their local analytic behavior away from such points.

Applying~\Cref{lem: maximum-modulus-principle} (Maximum Modulus Principle) to the holomorphic functions $\gamma_{\ell}(\lambda)$, we conclude that for each $\ell=1,\ldots,n$,
\begin{equation}
\label{eqn: gamma-k-infty-app-3}
|\gamma_{\ell}(\lambda)| \rightarrow \infty, \qquad \mathrm{as}\;\; |\lambda| \rightarrow \infty. 
\end{equation}
Furthermore, since $p'(z)$ is also a polynomial of degree $n-1$, it satisfies
\begin{equation}
\label{eqn: derivative-infty-app-3}
|p'(z)| \rightarrow \infty, \qquad \mathrm{as}\;\; |z| \rightarrow \infty.
\end{equation}
Combining~\eqref{eqn: gamma-k-infty-app-3} and~\eqref{eqn: derivative-infty-app-3}, we conclude that $|p'(\gamma_{\ell}(\lambda))| \rightarrow \infty$ as $|\lambda| \rightarrow \infty$. This establishes the uniform lower bound~\eqref{eqn: derivative-bound} and the asymptotic behavior~\eqref{eqn: derivative-bound-infty}, thereby completing the first part of~\Cref{lem: roots-seperation}.

\subsection*{Step 2: Lower bound on root seperation (compact region)}
Let $M>0$ be a constant to be determined later. Since the roots $\gamma_{\ell}(\lambda)$ vary continuously with $\lambda$, and the minimum separation function, 
\[
\min_{\ell_1 \neq \ell_2} | \gamma_{\ell_1}(\lambda) - \gamma_{\ell_2}(\lambda)| ,
\] 
is also continuous in $\lambda$. Therefore, by the compactness of the set $\overline{O(0,M)} \setminus S_p(\epsilon)$,  it follows that there exists a constant $\kappa(\epsilon) > 0$ such that 
\begin{equation}
\label{eqn: local-compact-app-3}
\min_{\ell_1 \neq \ell_2} | \gamma_{\ell_1}(\lambda) - \gamma_{\ell_2}(\lambda)| \geq \kappa(\epsilon), \quad \forall \lambda \in \overline{O(0,M)} \setminus S_p(\epsilon). 
\end{equation}

\subsection*{Step 3: Asymptotic estimates for large $|\lambda|$}
Since $p(z)$ is a degree-$n$ polynomial, we can write
\[
p(z) = a_nz^n \left( 1 + \frac{a_{n-1}}{a_nz} + \cdots + \frac{a_0}{a_nz^n}\right). 
\]
Hence, there exists a sufficiently large $R>0$ such that for any $|z| > R$, the lower-order terms are small, and we have the approximation
\[
\left| p(z) - a_nz^n \right| \leq \epsilon_0 |a_nz^n|,
\]
for some $\epsilon_0 \in (0,1/4)$. This implies the two-sided inequality:
\[
(1-\epsilon_0)|a_n| |z|^n \leq |p(z)| = |\lambda| \leq (1+\epsilon_0)|a_n| |z|^n.
\]
Solving for $|z|$ in terms of $|\lambda|$, we obtain bounds on the modulus of each root $\gamma_{\ell}(\lambda)$, for $\ell = 1,\ldots,n$:
\begin{equation}
\label{eqn: asymptotic-bound-app-3}
 \left( \frac{|\lambda|}{(1+\epsilon_0)|a_n|} \right)^{\frac1n} \leq |\gamma_{\ell}(\lambda)| \leq \left( \frac{|\lambda|}{(1-\epsilon_0)|a_n|} \right)^{\frac1n}.
\end{equation}

\subsection*{Step 4: Lower bound on root differences (asymptotic regime)}
The constant $M>0$ is chosen sufficiently such that the asymptotic bound~\eqref{eqn: asymptotic-bound-app-3} holds for all $|\lambda| \geq M$. In this asymptotic region,  we estimate the difference $\gamma_{\ell_1}(\lambda) - \gamma_{\ell_2}(\lambda)$ for any $\ell_1 \neq \ell_2$ as
\begin{align}
|\gamma_{\ell_1}(\lambda) - \gamma_{\ell_2}(\lambda) | & = \left|  \left( \frac{|\lambda|}{(1-\epsilon_0)|a_n|} \right)^{\frac1n} e^{\frac{i2\ell_1\pi}{n}} -  \left( \frac{|\lambda|}{(1+\epsilon_0)|a_n|} \right)^{\frac1n}e^{\frac{i2\ell_2\pi}{n}}  \right|   \nonumber \\
                                                                            & \geq  \left( \frac{|\lambda|}{(1-\epsilon_0)|a_n|} \right)^{\frac1n} \left| e^{\frac{i2\ell_1\pi}{n}} - e^{\frac{i2\ell_2\pi}{n}}   \right| - \left|  \left( \frac{|\lambda|}{(1-\epsilon_0)|a_n|} \right)^{\frac1n}  -  \left( \frac{|\lambda|}{(1+\epsilon_0)|a_n|} \right)^{\frac1n} \right|  \nonumber \\
                                                                            & \geq 2\left( \frac{|\lambda|}{(1-\epsilon_0)|a_n|} \right)^{\frac1n} \sin\left( \frac{(\ell_1 - \ell_2)\pi}{n}\right) - \frac{8\epsilon_0}{n}  \left| \frac{\lambda}{a_n} \right|^{\frac1n} \geq \frac{2}{n}  \left| \frac{\lambda}{a_n} \right|^{\frac1n}, \label{eqn: lower-bound-asymptotic-app-3}
\end{align}
where the first inequality follows the inequality for any $\epsilon_0 \in (0,1/4)$
\[
\left(\frac{1}{1-\epsilon_0}\right)^{\frac1n} - \left(\frac{1}{1+\epsilon_0}\right)^{\frac1n}  \leq \left(1 + 2\epsilon_0\right)^{\frac1n} - \left(1-2\epsilon_0\right)^{\frac1n}  \leq \frac{8\epsilon_0}{n}
\]
and the second inequality follows the identity 
\[
\left| e^{\frac{i2\ell_1\pi}{n}} - e^{\frac{i2\ell_2\pi}{n}}   \right|  = 2 \left| \sin \left( \frac{(\ell_1 - \ell_2)\pi}{n} \right) \right| \geq \frac{4}{n}.
\]
Combining this asymptotic lower bound~\eqref{eqn: lower-bound-asymptotic-app-3} with the compact-domain bound from~\eqref{eqn: local-compact-app-3}, we conclude that there exists $\kappa(\epsilon)>0$ such that  
\[
\min_{\ell_1 \neq \ell_2} | \gamma_{\ell_1}(\lambda) - \gamma_{\ell_2}(\lambda)|  \geq \kappa(\epsilon), \qquad \forall \lambda \in \mathbb{C} \setminus S_p(\epsilon), 
\]
which establishes the root-separation lower bound stated in~\eqref{eqn: root-seperation}, as well as the asymptotic behavior in~\eqref{eqn: root-seperation-infty}. This completes the proof of the second part of~\Cref{lem: roots-seperation}. 

\end{document}